# SPINAL PARTITIONS AND INVARIANCE UNDER RE-ROOTING OF CONTINUUM RANDOM TREES

By Bénédicte Haas, Jim Pitman[1] and Matthias Winkel[1,2]

*Université de Paris Dauphine, University of California and University of Oxford*

We develop some theory of spinal decompositions of discrete and continuous fragmentation trees. Specifically, we consider a coarse and a fine spinal integer partition derived from spinal tree decompositions. We prove that for a two-parameter Poisson–Dirichlet family of continuous fragmentation trees, including the stable trees of Duquesne and Le Gall, the fine partition is obtained from the coarse one by shattering each of its parts independently, according to the same law. As a second application of spinal decompositions, we prove that among the continuous fragmentation trees, stable trees are the only ones whose distribution is invariant under uniform re-rooting.

**1. Introduction.** Starting from a rooted combinatorial tree $T_{[n]}$ with $n$ leaves labeled by $[n] = \{1, \ldots, n\}$, we call the path from the root to the leaf labeled 1 the *spine* of $T_{[n]}$. Deleting each edge along the spine of $T_{[n]}$ defines a graph whose connected components we call *bushes*. If, as well as cutting each edge on the spine, we cut each edge connected to a spinal vertex, each bush is further decomposed into *subtrees*. We thus obtain two nested partitions of $\{2, \ldots, n\}$, which naturally extend to partitions of $[n]$ by adding the singleton $\{1\}$. We call these partitions of $[n]$ the *coarse spinal partition* and the *fine spinal partition* derived from $T_{[n]}$. See, for example, Figure 2.

The aim of this paper is to develop some theory of spinal decompositions of *fragmentation trees* that arise as genealogical trees of fragmentation processes. We focus on Markovian partition-valued fragmentation processes of the following two types. In a setting of discrete time and partitions of $[n]$,

Received May 2007; revised March 2008.
[1]Supported in part by EPSRC Grant GR/T26368/01 and NSF Grant DMS-04-05779.
[2]Supported by the Institute of Actuaries and Aon Limited.

we postulate that each nonsingleton block splits at each time, which leads to Markov branching models [4, 18, 26]. In a setting of continuous time and partitions of $\mathbb{N}$ we postulate a self-similarity condition, which leads to self-similar continuum random trees [25, 26].

Before giving an overview of this paper in Section 1.3, we formally introduce the discrete setting in Section 1.1 and the continuous setting in Section 1.2.

1.1. *Discrete fragmentations.* We start by introducing a convenient formalism for the kind of combinatorial trees arising naturally in the context of fragmentation processes. Let $B$ be a finite nonempty set, and write $\#B$ for the number of elements of $B$. Following standard terminology, a *partition of* $B$ is a collection

$$\Pi_B = \{B_1, \ldots, B_k\}$$

of nonempty disjoint subsets of $B$ whose union is $B$. To introduce a new terminology convenient for our purpose, we make the following recursive definition. A *fragmentation of* $B$ (sometimes called a *hierarchy* or a *total partition* [35, 36]) is a collection $T_B$ of nonempty subsets of $B$ such that:

   (i) $B \in T_B$,
   (ii) if $\#B \geq 2$ there is a partition $\Pi_B$ of $B$ into at least two parts $B_1, \ldots, B_k$, called the *children of* $B$, with

(1)  $$T_B = \{B\} \cup T_{B_1} \cup \cdots \cup T_{B_k},$$

where $T_{B_i}$ is a fragmentation of $B_i$ for each $1 \leq i \leq k$.

Necessarily $B_i \in T_B$, each child $B_i$ of $B$ with $\#B_i > 1$ has further children, and so on, until the set $B$ is broken down into singletons. We use the same notation $T_B$ for both:

- such a collection of subsets of $B$, and
- for the tree whose vertices are these subsets of $B$, and whose edges are defined by the parent/child relation implicitly determined by the collection of subsets of $B$.

To emphasize the tree structure we may call $T_B$ a *fragmentation tree*. Thus $B$ is the first branch point of $T_B$, and each singleton subset of $B$ is a leaf of $T_B$, see Figure 1. It is often convenient to plant $T_B$ by adding a ROOT vertex and an edge between the ROOT and the first branch point $B$. We denote by $\mathbb{T}_B$ the collection of all fragmentation trees labeled by $B$.

For each nonempty subset $A$ of $B$, the *restriction to $A$ of $T_B$*, denoted $T_{A,B}$, is the fragmentation tree whose first branch point is $A$, whose leaves are the singleton subsets of $A$, and whose tree structure is defined by restriction of $T_B$. That is, $T_{A,B}$ is the fragmentation

$$T_{A,B} = \{C \cap A : C \cap A \neq \varnothing, C \in T_B\} \in \mathbb{T}_A,$$

corresponding to a *reduced subtree* as discussed by Aldous [1].

Given a rooted combinatorial tree with no single-child vertices and whose leaves are labeled by a finite set $B$, there is a corresponding fragmentation tree $T_B$, where each vertex of the combinatorial tree is associated with the set of leaves in the subtree above that vertex. So the fragmentation trees defined here provide a convenient way to both label the vertices of a combinatorial tree, and to encode the tree structure in the labeling.

A *random fragmentation model* is an assignment of a probability distribution on $\mathbb{T}_B$ for a random fragmentation tree $T_B$ with first branch point $B$ for each finite subset $B$ of $\mathbb{N}$. We assume throughout this paper that the model is *exchangeable*, meaning that the distribution of $\Pi_B$, the partition of

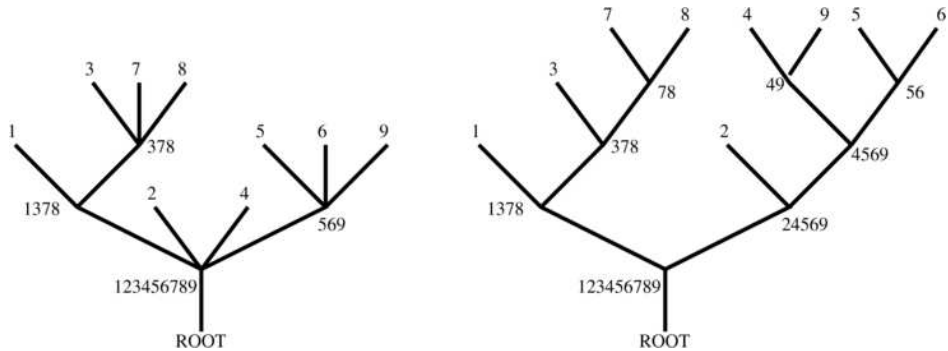

Fig. 1. *Two fragmentations of [9] represented as trees with nodes labeled by subsets of [9].*

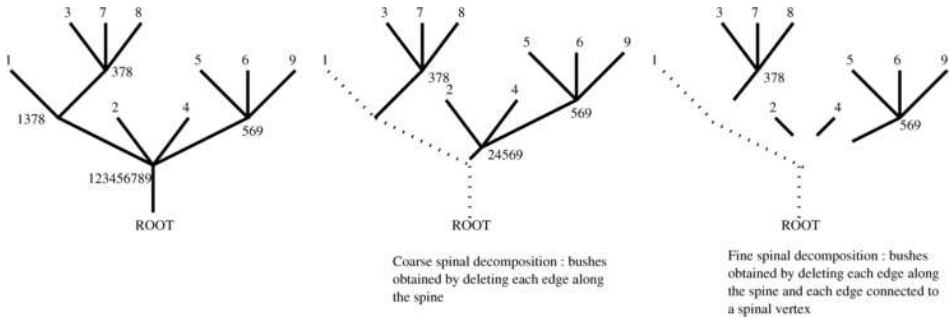

Fig. 2. *A fragmentation tree $T_{[9]}$, with coarse spinal partition $\{\{1\}, \{24569\}, \{378\}\}$, coarse spinal composition $(\{24569\}, \{378\}, \{1\})$ and fine spinal partition $\{\{1\}, \{2\}, \{378\}, \{4\}, \{569\}\}$.*

$B$ generated by the branching of $T_B$ at its root, is of the form

(2) $$\mathbb{P}(\Pi_B = \{B_1, \ldots, B_k\}) = p(\#B_1, \ldots, \#B_k)$$

for all partitions $\{B_1, \ldots, B_k\}$ with $k \geq 2$ blocks, and some symmetric function $p$ of compositions of positive integers, called a *splitting probability rule*. The model is called:

- *Markovian* (or a Markov branching model) if given $\Pi_B = \{B_1, \ldots, B_k\}$, the $k$ subtrees of $T_B$ above $B$ are independent and distributed as $T_{B_1}, \ldots, T_{B_k}$, for all partitions $\{B_1, \ldots, B_k\}$ of $B$;
- *consistent* if for every $A \subset B$, the restriction to $A$ of $T_B$ is distributed like $T_A$;
- *binary* if every $A \in T_B$ has either 0 or 2 children with probability one, for all $B$.

Now we take $B = [n]$. The collection of vertices at graph distance $m \geq 0$ above the first branch point form a partition of a subset of $[n]$ that we extend to a partition $\Pi_m^{(n)}$ of $[n]$ by adding a singleton $\{j\}$ for each leaf $j$ at height below $m$. We refer to $(\Pi_m^{(n)}, m \geq 0)$ as the partition-valued discrete *fragmentation process* associated with $T_{[n]}$. See also [4, 18, 26].

1.2. *Continuous self-similar fragmentations.* We denote by $\mathcal{P}$ the set of partitions of $\mathbb{N}$ and equip it with the distance $d(\pi, \pi') = 2^{-n(\pi, \pi')}$, where $n(\pi, \pi')$ is the largest integer such that the restrictions of partitions $\pi, \pi'$ to $[n]$ coincide. Following Bertoin [9], a continuous-time $\mathcal{P}$-valued Markov process $(\Pi(t), t \geq 0)$ is called a *self-similar fragmentation process* with index $a \in \mathbb{R}$ if it is *càdlàg* and:

- $\Pi(0) = \{\mathbb{N}\}$, that is, $\Pi$ starts from the trivial partition with a unique block;
- $\Pi$ is exchangeable, that is, its distribution is invariant under permutations of $\mathbb{N}$;
- given $\Pi(t) = \pi$, the post-$t$ process $(\Pi(t+s), s \geq 0)$ has the same law as the process whose state at time $s \geq 0$ is the partition of $\mathbb{N}$ whose blocks are those of

$$\pi_i \cap \Pi^{(i)}(|\pi_i|^a s), \qquad i \geq 1,$$

where $(\pi_i, i \geq 1)$ is the sequence of blocks of $\pi$ in order of least elements, $(|\pi_i|, i \geq 1)$ is the sequence of their asymptotic frequencies and $(\Pi^{(i)}, i \geq 1)$ is a sequence of i.i.d. copies of $\Pi$.

We recall that Kingman's theorem [27] on exchangeable partitions ensures that for every $t \geq 0$, the asymptotic frequencies $|\pi_i| = \lim_{n \to \infty} n^{-1} \#(\pi_i \cap [n])$ of all blocks $\pi_i$ of $\Pi(t)$ exist a.s. Bertoin [8] shows that actually a.s. for every $t$, these asymptotic frequencies exist.

In [9], Bertoin proved that the distribution of a self-similar fragmentation is entirely characterized by three parameters: the *index of self-similarity* $a$, a coefficient $c \geq 0$ that measures the rate of *erosion* and a *dislocation measure* on

$$\mathcal{S}^{\downarrow} = \left\{ (s_i)_{i \geq 1} : s_1 \geq s_2 \geq \cdots \geq 0, \sum_{i \geq 1} s_i \leq 1 \right\}$$

with no atom at $(1, 0, \ldots)$ and that integrates $1 - s_1$. This measure $\nu$ describes the sudden dislocations of blocks, in the sense that a block $B \subset \mathbb{N}$ splits in some blocks $B_1, B_2, \ldots$ with relative asymptotic frequencies $\mathbf{s} \in \mathcal{S}^{\downarrow}$ at rate $|B|^a \nu(d\mathbf{s})$. When the index $a = 0$, this fragmentation rate does not depend on the size of the blocks and the fragmentation processes is then said to be *homogeneous*. A crucial point is that a self-similar fragmentation with parameters $a$, $c$ and $\nu$ can always be constructed measurably from a homogeneous fragmentation with same coefficient $c$ and measure $\nu$, using time-changes, and vice versa. We refer to Bertoin's book [10] and the above mentioned papers [8, 9] for details on these time-changes and background on homogeneous and self-similar fragmentations.

In this paper, we focus on self-similar fragmentations without erosion ($c = 0$), which are nontrivial ($\nu(\mathcal{S}^{\downarrow}) \neq 0$) and do not lose mass at sudden dislocations, that is,

$$\nu\left( \sum_{i \geq 1} s_i < 1 \right) = 0.$$

We call $(a, \nu)$ the characteristic pair of such a process.

A family of combinatorial trees with edge lengths $R_{[n]}, n \geq 1$, with $n$ exchangeably labeled leaves, is naturally associated to a self-similar fragmentation process $\Pi$ by considering the evolution of $\Pi$ restricted to the first $n$ integers. Specifically, $R_{[n]}$ consists of all blocks $B$ that occur in the evolution of $\Pi \cap [n]$; an edge between the root and the first branch point $[n]$ has as its length the first dislocation time of $\Pi \cap [n]$, and similarly for subtrees with two or more leaves; the edge below leaf $j$ has as its length the time between the last relevant dislocation time of $\Pi \cap [n]$ and the time when $\{j\}$ becomes a singleton for $\Pi$, which may be infinite. This gives a *consistent* family of trees, in the sense that the subtree of $R_{[n]}$ spanned by $[k]$ is $R_{[k]}$, for all $k \leq n$, where superfluous (i.e., multiplicity 2) vertices are removed and associated edges merged, their lengths summed up. By exchangeability, the same is true in distribution for uniformly chosen $k$ distinct leaves of $R_{[n]}$, relabeled by $[k]$. The coupling of self-similar fragmentations using time-changes entails that the distribution of the combinatorial shapes (say $T_{[n]}$) of $R_{[n]}, n \geq 1$, depends only on the dislocation measure $\nu$, and not on the index $a$. So without loss of generality, we may focus on $a = 0$, the case of homogeneous fragmentations, when working with the shapes $T_{[n]}, n \geq 1$. Furthermore, $(T_{[n]}, n \geq 1)$ defines a consistent Markov branching model as in the previous subsection. Reciprocally, each consistent Markov branching model can be constructed similarly from some homogeneous fragmentation (possibly with erosion). See [26].

When the index $a$ is negative, small fragments vanish quickly and it is well known that the whole fragmentation $\Pi$ then reaches in finite time the trivial partition composed exclusively of singletons. See, for example, [10]. In terms of trees, this implies that the height of $R_{[n]}$ is bounded uniformly in $n$. Using the consistency property and Aldous' results [3], it is then possible to define the projective limit $\mathcal{T}$ of the family $(R_{[n]}, n \geq 1)$ and equip it with a probability measure $\mu$, the *mass measure*, that arises as limit of the empirical measures on the leaves of $R_{[n]}, n \geq 1$. Implicitly, the tree $\mathcal{T}$ is rooted. The pair $(\mathcal{T}, \mu)$ is a *continuum random tree* (CRT) and was studied in [25] using Aldous' formalism of trees as compact metric subsets of $l_1$; cf. [1, 2, 3]. An alternative formalism can be considered, via the set of equivalence classes of compact rooted $\mathbb{R}$-trees endowed with the Gromov–Hausdorff distance, as developed in [16, 17]. We will not go further into details here and refer to the above-mentioned papers for rigorous definitions and statements. We shall call the CRT $(\mathcal{T}, \mu)$ a *self-similar fragmentation CRT* with parameters $(a, \nu)$.

A fundamental property of $(\mathcal{T}, \mu)$ is that a version of $(R_{[n]}, n \geq 1)$ can be obtained from a random sampling $L_1, L_2, \ldots$ picked independently according to $\mu$, conditional on $(\mathcal{T}, \mu)$, by considering for each $n$ the subtree of $\mathcal{T}$ spanned by the root and leaves $L_1, \ldots, L_n$. Consider then the forest $\mathcal{F}_\mathcal{T}(t)$ obtained by removing in $\mathcal{T}$ all vertices at distance less than $t$ from the root and define the random partition $\Pi'(t)$ by letting $i$ and $j$ be in the same block of $\Pi'(t)$ if and only if $L_i$ and $L_j$ are in the same connected component of $\mathcal{F}_\mathcal{T}(t)$, $t \geq 0$. Clearly the process $\Pi'$ is distributed as $\Pi$. We shall often suppose in the following that the fragmentation process we are working with is constructed in such a manner from some self-similar fragmentation CRT.

Examples of self-similar fragmentation CRTs are the Brownian CRT of Aldous [1, 2, 3] and, more generally, the stable Lévy trees with index $\beta \in (1, 2]$ of Duquesne and Le Gall [14, 15]. For details on their fragmentation properties, see Bertoin [9] for the Brownian case (i.e., when $\beta = 2$) and Miermont [29] for the other stable cases. The parameters of these CRTs are recalled later in the paper.

1.3. *Contents and organization of the paper.* The structure and contents of this paper are as follows. In Section 2, we study the coarse and fine spinal partitions of some Markov branching model $(T_{[n]}, n \geq 1)$ constructed consistently from a self-similar fragmentation process. These partitions of $[n]$ are consistent as $n$ varies, which leads to a nested pair of partitions of $\mathbb{N}$. Restricted to $\mathbb{N} \setminus \{1\}$, they are jointly exchangeable. In particular, they possess asymptotic frequencies a.s. The decreasing rearrangements of these frequencies are called the *coarse spinal mass partitions* and *fine spinal mass partitions*. By decomposing the trees along the spine, we then show that when the parameters $a$ and $\nu$ of the fragmentation are known and $\nu$ is infinite, we can reconstruct the whole self-similar fragmentation process from the sequence of shapes $(T_{[n]}, n \geq 1)$ (Proposition 2). Next, the main result of this section (Theorem 6) states that under some factorization property of the dislocation measure $\nu$ (Definition 2), the fine spinal mass partition derived from the sequence of shapes $(T_{[n]}, n \geq 1)$ is obtained from the coarse one by shattering each of its fragments in an i.i.d. manner.

In particular, this result applies to a family of fragmentations whose dislocation measures are built from Poisson–Dirichlet partitions (Section 3). The stable fragmentations studied by Miermont [29], built from the stable trees of Duquesne and Le Gall with index in $(1, 2)$, belong to this family. As a consequence, we obtain an extensive description, in terms of Poisson–Dirichlet partitions (Corollary 10), of spinal decompositions of stable trees.

The stable trees $(\mathcal{T}, \mu)$ are known to possess an interesting symmetry property of *invariance under uniform re-rooting* (see [2, 13, 15]). Informally, this means that taking a leaf at random according to $\mu$ and considering $\mathcal{T}$ rooted at this random leaf, gives a CRT with the same distribution as the original CRT with its original root. In Section 4, we give a new proof of this invariance, using combinatorial methods, and show that, up to a scaling factor, stable trees are the only self-similar fragmentation CRTs that are invariant under uniform re-rooting (Theorem 11).

To finish this introduction, let us mention that studies on spinal decompositions of various trees exist in the literature. See, for example, Aldous–Pitman [6] (for Galton–Watson trees), Duquesne–Le Gall [15] (for stable and Lévy trees). In the fragmentation context, decomposing the trees/processes along the spine is a useful tool, which has been used to obtain results on large time asymptotics [11], small time asymptotics [24] and discrete approximations [26].

**2. Spinal partitions of fragmentation trees.** Decompose a combinatorial fragmentation tree $T_{[n]}$ with leaves labeled by $[n]$ along the spine from the root to leaf 1 into a collection of bushes by deleting each edge along the spine. By adding a conventional root edge to its base, each bush is identified with an element of $\mathbb{T}_B$ for some $B \subseteq [n]$, where $\mathbb{T}_B$ is the collection of rooted combinatorial trees with $\#B$ leaves labeled by $B$. Each such $B$ is associated with a unique vertex on the spine of $T_{[n]}$. We list these sets of leaf labels $B$ in order of the corresponding spinal vertices to obtain an *ordered* exchangeable random partition of $\{2, \ldots, n\}$. The first set $B$ in this list is the set of elements of $[n]$ not in the block containing 1 after the first fragmentation event involving $[n]$. If after the first fragmentation of $[n]$ the block $[n] - B$ containing 1 is of size 2 or more, the next set is what remains of $[n] - B$ after deleting the block containing 1 in the next fragmentation of $[n] - B$, and so on, until the last set which is the singleton $\{1\}$. If as well as cutting each edge on the spine, we cut each edge connected to a spinal vertex, each bush is further decomposed into subtrees. We thus obtain two nested exchangeable random partitions of $\{2, \ldots, n\}$, which naturally extend to partitions of $[n]$ by adding the singleton $\{1\}$, the coarse and fine spinal partitions derived from $T_{[n]}$. We can include the spinal order in the coarse spinal partition to form the *coarse spinal composition*.

Assuming that the trees $T_{[n]}, n \geq 1$, are constructed consistently from a homogeneous fragmentation process with values in the partitions of $\mathbb{N}$, both partitions of $[n]$ are consistent as $n$ varies. Thus the coarse and fine spinal partitions may be regarded as a nested pair of random partitions of $\mathbb{N}$. These partitions have natural interpretations in terms of associated partition-valued *self-similar* fragmentations processes $(\Pi(t), t \geq 0)$, of any index $a$, in which the sequence $(T_{[n]}, n \geq 1)$ is embedded. For each pair of integers $i$ and $j$ let the *splitting time* $D_{i,j}$ be the first time $t$ that $i$ and $j$ fall in distinct blocks of $\Pi(t)$. Let $i, j \geq 2$. By construction, $i$ and $j$ fall in the same block of the coarse spinal partition if and only if $D_{1,i} = D_{1,j}$, whereas $i$ and $j$ fall in the same block of the fine spinal partition if and only if $D_{i,j} > D_{1,i}$ (this clearly implies $D_{1,i} = D_{1,j}$). Assuming further that $\Pi$ is constructed by random sampling of leaves $L_1, L_2, \ldots$ from some CRT $(\mathcal{T}, \mu)$ according to $\mu$, $i$ and $j$ fall in the same block of the coarse spinal partition if and only if the paths from $L_i$ and $L_j$ to the root first meet the *spine of $\mathcal{T}$*, that is, the path from the root to $L_1$, at the same point. Besides, $i$ and $j$ fall in the same block of the fine spinal partition if and only if the path from $L_i$ to $L_j$ does not intersect the spine.

The *coarse spinal decomposition of $\mathcal{T}$* is the collection of equivalence classes for the random equivalence relation $x \sim y$ if and only if the paths from $x$ and $y$ to the root first meet the spine at the same point on the spine. Note that the whole spine itself carries no $\mu$-mass, and spinal non-branchpoints (an uncountable set of singletons in this decomposition of $\mathcal{T}$) will be excluded from further consideration. The restriction of $\mathcal{T}$ to a typical equivalence class is a bush which can be further decomposed into trees by deleting the point on the spine, and then giving each connected component its own root where it used to be connected to the spine. The resulting random partition of $\mathcal{T}$ into subtrees is the *fine spinal decomposition of $\mathcal{T}$*.

We measure the size of each component of one of these partitions by its $\mu$-mass, to obtain *coarse and fine spinal mass partitions of $(\mathcal{T}, \mu)$*, which we may regard as two random elements of $\mathcal{S}^{\downarrow}$. The following proposition summarizes some basic properties of these random partitions, which follow easily from the above discussion.

PROPOSITION 1. *The coarse and fine spinal partitions derived from the sequence of shapes $(T_{[n]}, n \geq 1)$ embedded in $(\mathcal{T}, \mu)$ have the following properties:*

(i) *The singleton block $\{1\}$ belongs to both partitions of $\mathbb{N}$, while the restrictions of these partitions to $\mathbb{N} \setminus \{1\}$ are jointly exchangeable.*

(ii) *The sequence of ranked limiting frequencies of each partition of $\mathbb{N}$ is the sequence of ranked $\mu$-masses of the corresponding mass partition of $(\mathcal{T}, \mu)$.*

We now offer a more detailed study of these two partitions, first considering the coarse spinal partition (and composition), then the fine one and its relation to the coarse one. Obviously, the fine spinal partition is identical to the coarse one if and only if the trees $T_{[n]}$ are binary for all $n \geq 1$.

2.1. *The coarse spinal partition.* Assume throughout this section that the trees $T_{[n]}, n \geq 1$, are constructed consistently from a homogeneous fragmentation process, as when $T_{[n]}$ is derived from an associated continuum tree $(\mathcal{T}, \mu)$ as the shape of the subtree spanned by $L_i, i \in [n]$, for $L_1, L_2, \ldots$ an exchangeable sample of leaves with directing measure $\mu$. To ease notation we work with $T_{[n+1]}$ instead of $T_{[n]}$. Let

$$B_{n,1}, B_{n,2}, \ldots, B_{n,K_n}, \{1\}$$

be the sets of leaves of the bushes derived from the coarse spinal decomposition of $T_{[n+1]}$, in order of the corresponding spinal vertices. Then $(B_{n,1}, B_{n,2}, \ldots, B_{n,K_n})$ is the restriction to $\{2, \ldots, n+1\}$ of an exchangeable ordered random partition of $\{2, 3, \ldots\}$, as studied in [12, 21]. Let

(3) $$\mathcal{C}_n := (\#B_{n,1}, \#B_{n,2}, \ldots, \#B_{n,K_n}).$$

It follows easily from sampling consistency of the sequence $(T_{[n]}, n \geq 1)$ that $(\mathcal{C}_n, n \geq 1)$ is a *regenerative composition structure*, as defined in [19]. That is to say, $(\mathcal{C}_n, n \geq 1)$ is a sampling consistent sequence of random compositions $\mathcal{C}_n$ of $n$, with the property that conditionally given the first part of $\mathcal{C}_n$ is of size $i < n$, the remaining parts of $\mathcal{C}_n$ define a random composition of $n - i$ with the same distribution as $\mathcal{C}_{n-i}$. Let

$$S_{n,k} := n - \sum_{j=1}^{k} \#B_{n,j},$$

where $B_{n,j}$ is empty for $j > K_n$. So $(S_{n,k} + 1, k \geq 0)$ is the sequence of sizes of the fragment containing 1 as it undergoes successive fragmentations according to $T_{[n+1]}$, starting with $S_{n,0} = n$ and terminating with $S_{n,k} = 0$ for $k \geq K_n$, where $K_n$ is the total number of fragmentation events experienced by the block containing 1 in $T_{[n+1]}$. According to Gnedin and Pitman [19], there is the following almost sure convergence of random sets with respect to the Hausdorff metric on closed subsets of $[0, 1]$:

$$
(4) \qquad \{S_{n,k}/n, k \geq 0\} \xrightarrow[n \to \infty]{\text{a.s.}} \{\exp(-\xi_t), t \geq 0\}^{\text{cl}},
$$

where the left-hand side is the random discrete set of values $S_{n,k}$ rescaled onto $[0, 1]$, and the right-hand side is the closure of the range of the exponential of some subordinator $(\xi_t, t \geq 0)$. The random interval partition of $[0, 1]$ defined by interval components of the complement of the closed range of $1 - e^{-\xi}$ has a natural interpretation in terms of the associated CRT $(\mathcal{T}, \mu)$: the lengths of these intervals are the strictly positive masses of components in the coarse spinal decomposition of $(\mathcal{T}, \mu)$, in the order they appear along the spine from the root to leaf 1. We will therefore call this interval partition the *coarse spinal interval partition of* $[0, 1]$ *derived from* $(\mathcal{T}, \mu)$. In terms of the associated homogeneous fragmentation, the lengths of these intervals are the total masses of fragments thrown off by the mass process of the fragment containing 1, put in the order they split away from this tagged fragment. Otherwise said, $\exp(-\xi)$ is the mass process of the fragment containing 1. Since the fragmentation process has zero erosion and no sudden loss of mass, the subordinator $\xi$ has no drift and no killing. Bertoin [8] proved that the Lévy measure of $\xi$ is then given by

$$
(5) \qquad \Lambda(dx) = \exp(-x) \sum_{i \geq 1} \nu(-\log s_i \in dx), \qquad x > 0.
$$

PROPOSITION 2. *Let* $(\Pi(t), t \geq 0)$ *be a self-similar fragmentation process, with index* $a \in \mathbb{R}$ *and dislocation measure* $\nu$ *with infinite total mass. Then the entire process* $(\Pi(t), t \geq 0)$ *can be constructed from the consistent sequence* $(T_{[n]}, n \geq 1)$ *of combinatorial shapes of trees derived from* $\Pi$.

PROOF. In view of the time-change relation between fragmentations of different indices, it suffices to consider the homogeneous case. Given the consistent family of trees $(T_{[n]}, n \geq 1)$, we first use (4) to recover the closure of the range of $\exp(-\xi)$, hence also the closure of the range of $\xi$, the subordinator describing the evolution of the mass fragment containing 1. Since the dislocation measure has infinite mass, so does the Lévy measure of $\xi$. Then it is well known that the entire sample path of $\xi$ can be measurably reconstructed from its range, up to a constant factor on the time scale (see, e.g., [22]). Since the distribution of $\xi$ is determined by that of $(\Pi(t), t \geq 0)$, this constant is known. Let $\Pi_n = (\Pi_n(t), t \geq 0)$ be the restriction of $(\Pi(t), t \geq 0)$ to $[n]$. The path of $\xi$, and its construction (4) from $(T_{[n]}, n \geq 1)$, determine almost surely for each $n$ the sequence of random times $t$ when transitions of $\Pi_n$ occur which change the block of $\Pi_n$ containing 1, and at each of these times $t$ the block of $\Pi_n(t)$ containing 1 can be read from $T_{[n]}$. By exchangeability, the same reconstruction can evidently be done almost surely for the block of $\Pi_n(t)$ containing $j$, for each $1 \leq j \leq n$. But this information determines the entire path of $(\Pi_n(t), t \geq 0)$, for each $n$, hence that of $(\Pi(t), t \geq 0)$. □

COROLLARY 3. *If in the setting of Proposition 2 we have $a < 0$, then an associated $(a, \nu)$-fragmentation CRT $(\mathcal{T}, \mu)$ can also be constructed from $(T_{[n]}, n \geq 1)$ on the same probability space.*

PROOF. While the construction of a self-similar fragmentation CRT in [25] from a self-similar partition-valued fragmentation process is carried out explicitly only "in distribution," it is not hard to give an almost sure construction, for example, via Aldous' sequential construction in $l_1$ (see, e.g., [3], page 252). This yields an increasing sequence of trees with edge lengths $R_{[n]}$ that converges in distribution, hence almost surely, with respect to the Hausdorff metric on closed subsets of $l_1$. The almost sure convergence of empirical measures on the leaves of $R_{[n]}$ to a mass measure $\mu$ is then given by [3], Lemma 7 (convergence of measures is weak convergence). □

We record now an explicit distributional result for the coarse spinal partition of $T_{[n+1]}$, which can either be read from [19] or derived directly. Recall that $n + 1 - \#B_{n,1}$ is the size of the fragment containing 1 at the first branch point of $T_{[n+1]}$. Let $\Sigma(ds) := \sum_{j=1}^{\infty} \nu(s_j \in ds)$ and let $\Lambda$ be the Lévy measure of $(\xi_t, t \geq 0)$, which, according to (5), is the image of $s\Sigma(ds)$ via $s \mapsto -\log s$. Then by embedding in the homogeneous fragmentation, we see that

$$\mathbb{P}(\#B_{n,1} = m) = \Phi(n:m)/\Phi(n) \qquad (1 \leq m \leq n), \tag{6}$$

where $\Phi(n)$ is the total rate of fragmentations with some effect on partitions of $[n+1]$, and $\Phi(n:m)$ the rate of such fragmentations in which 1 ends up in a block of size $n + 1 - m$. From standard results on the construction of the homogeneous fragmentation from its dislocation measure $\nu$ (cf. [10], Chapter 3), these rates are easily evaluated as follows:

$$\Phi(n:m) = \binom{n}{m} \int_0^1 s^{n+1-m}(1-s)^m \Sigma(ds)$$
$$= \binom{n}{m} \int_0^\infty e^{-(n-m)x}(1-e^{-x})^m \Lambda(dx) \tag{7}$$

and

$$\Phi(n) = \sum_{m=1}^n \Phi(n:m) = \int_0^1 (1 - s^n) s \Sigma(ds) = \int_0^\infty (1 - e^{-nx}) \Lambda(dx). \tag{8}$$

From this and [19], we get the exchangeable partition probability function (EPPF) of the coarse spinal partition $\{B_{n,1}, B_{n,2}, \ldots, B_{n,K_n}\}$ restricted to $\{2, \ldots, n+1\}$, that is, the probabilities $p(n_1, \ldots, n_k) = \mathbb{P}(\{B_{n,1}, \ldots, B_{n,K_n}\} = \pi)$, for each particular partition $\pi$ of $\{2, \ldots, n+1\}$ in sets of sizes $n_1, \ldots, n_k$, $\forall n \geq 1, \forall (n_1, \ldots, n_k)$ partition of $n$. For an explicit formula, see [19], especially formulae (3), (4), (6) and (26). Various further properties of the coarse spinal partition can also be read from [19].

2.2. *The fine spinal partition.* We start by observing some basic symmetry properties of this partition.

PROPOSITION 4. (i) *Consider the fine spinal partition derived from $T_{[n+1]}$, restricted to $\{2, \ldots, n+1\}$. Then, conditionally given the sizes of its components, say $n_1, \ldots, n_k$ with $\sum_{i=1}^{k} n_i = n$, the corresponding collection of subtrees of $T_{[n+1]}$, modulo relabeling by $[n_1], \ldots, [n_k]$, is a collection of independent copies of $T_{[n_1]}, \ldots, T_{[n_k]}$.*

(ii) *Conditionally given the fine spinal mass partition of a self-similar fragmentation CRT $(\mathcal{T}, \mu)$ with parameters $(a, \nu)$, the corresponding collection of subtrees $T$ of $\mathcal{T}$, with each $T$ of mass $m$ equipped with $m^{-1}\mu$ restricted to $T$, modulo isomorphism and multiplication of edge lengths by $m^a$, is a collection of independent copies of $(\mathcal{T}, \mu)$.*

PROOF. Part (i) follows easily from the defining Markov (fragmentation/branching) property of $T_{[n]}$. For part (ii), consider $\Pi$ a partition-valued $(a, \nu)$-fragmentation constructed from $(\mathcal{T}, \mu)$. Let $\Pi_{(i)}(t)$ denote the block of $\Pi(t)$ containing $i$, $i \geq 1$, and recall that $D_{1,i}$ denotes the first time at which 1 and $i$ belong to distinct blocks. For all $t \geq 0$, the collection of blocks $(\Pi_{(i)}(D_{1,i}+t), i \geq 1)$ induces a partition of $\mathbb{N}$. In the terminology of Bertoin ([10], Definition 3.4), the sequence $(D_{1,i}, i \geq 1)$ is a *stopping line*, and as such, satisfies the *extended branching property* ([10], Lemma 3.14), which ensures that given $(\Pi_{(i)}(D_{1,i}), i \geq 1)$, the processes $(\Pi_{(i)}(D_{1,i}+t), t \geq 0)$, $i \geq 1$, evolve, respectively, as $(m_i \Pi^{(i)}(m_i^a t), t \geq 0)$, where $m_i$ is the asymptotic frequency of $\Pi_{(i)}(D_{1,i}), i \geq 1$, and the $\Pi^{(i)}$s are i.i.d. copies of $\Pi$. Now, coming back to the CRT $(\mathcal{T}, \mu)$, each component of its fine spinal partition corresponds to a fragmentation $(\Pi_{(i)}(D_{1,i}+t), t \geq 0)$ for some $i$ and obviously, can be measurably reconstructed from this fragmentation (see the proof of Corollary 3). Conditionally given the masses $m_i, i \geq 1$, the subtrees of the fine spinal partition are therefore independent, distributed (modulo isomorphisms), respectively, as $(m_i^{-a}\mathcal{T}, m_i \mu^{(m_i^{-a})})$, $i \geq 1$, where $m_i^{-a}\mathcal{T}$ means that the edge lengths of $\mathcal{T}$ have been multiplied by $m_i^{-a}$ and $\mu^{(m_i^{-a})}$ is the image of $\mu$ by this transformation. □

Part (ii) of the proposition is a natural generalization of the spinal decomposition of the Brownian CRT described in [5]. When the Brownian CRT is encoded in a Brownian excursion, this corresponds to a path decomposition whereby a single excursion is decomposed into a countably infinite collection of independent copies of itself.

In view of this symmetry property of the fine spinal partition, it is natural to look for some more explicit description of this decomposition, such as its EPPF or the distribution on $\mathcal{S}^\downarrow$ of the corresponding mass partition. While such descriptions are known for the Brownian CRT, and more generally for all binary self-similar fragmentation CRTs according to the previous section, they appear to be difficult to obtain in general. But searching for conditions which simplify the structure of the fine spinal partition of $(\mathcal{T}, \mu)$ leads naturally to consideration of further symmetry properties, and then to interesting examples with these properties for which explicit computations can be made. Consider first the fine partition of the set of leaves in some block of the coarse spinal partition of $T_{[n+1]}$ (restricted to $\{2, \ldots, n+1\}$). By recursive arguments, it is enough to discuss the fine partition of the first block of the coarse spinal partition.

For each $\mathbf{s} \in \mathcal{S}^\downarrow$ let $\mathbb{P}_\mathbf{s}$ denote the probability measure governing an exchangeable random partition $\Pi$ of $\mathbb{N}$ whose ranked frequencies are equal to $\mathbf{s}$, and for a measure $\nu$ on $\mathcal{S}^\downarrow$ let

$$\mathbb{P}_\nu(\cdot) = \int_{\mathcal{S}^\downarrow} \mathbb{P}_\mathbf{s}(\cdot) \nu(d\mathbf{s})$$

be the corresponding distribution of $\Pi$ as a mixture of Kingman's paintbox partitions. For each $n$ the distribution of $\Pi_n$ is determined by the formula

$$\mathbb{P}_\nu(\Pi_n = \{B_1, \ldots, B_k\}) = p_\nu(\#B_1, \ldots, \#B_k)$$

for every partition $\{B_1, \ldots, B_k\}$ of $[n]$ into $k \geq 1$ parts, for some function $p_\nu(n_1, \ldots, n_k)$ of compositions $(n_1, \ldots, n_k)$ of positive integers $n$. We refer here to [32] or [10] for a specific formula for $p_\nu(n_1, \ldots, n_k)$. In particular, $p_\nu(1,1) = \int_{\mathcal{S}^\downarrow}(1 - \sum_{i \geq 1} s_i^2) \nu(d\mathbf{s})$. Note that $p_\nu(n_1, \ldots, n_k) < \infty$ for all $n_1, \ldots, n_k \in \mathbb{N}$, $k \geq 2$, if and only if $p_\nu(1,1) < \infty$, that is, if and only if $\int_{\mathcal{S}^\downarrow}(1 - s_1) \nu(d\mathbf{s}) < \infty$.

DEFINITION 1. The function $p_\nu$ is called the *exchangeable partition rate function* (*EPRF*) associated with $\nu$. If $\nu$ is a probability measure, then so is $\mathbb{P}_\nu$, and $p_\nu$ is known as an *exchangeable partition probability function* (*EPPF*).

Note that we have the addition rule

$$p_\nu(n_1, \ldots, n_k) = p_\nu(n_1 + 1, n_2, \ldots, n_k) + \cdots + p_\nu(n_1, \ldots, n_{k-1}, n_k + 1)$$
$$+ p_\nu(n_1, \ldots, n_k, 1).$$

The following lemma presents a basic decomposition in some generality.

LEMMA 5. *Let $\nu$ be a dislocation measure on $\mathcal{S}^\downarrow$ with associated EPRF $p_\nu$. Then for every $k \geq 2$ and every composition $n_1, \ldots, n_k$ of $n \geq 2$ into at least two parts,*

$$p_\nu(n_1, \ldots, n_k) = g(n, n_1) p_{\widehat{\nu}(n,n_1)}(n_2, \ldots, n_k) \tag{9}$$

*for some function $g(n, n_1)$ and some family of probability measures $\widehat{\nu}(n, n_1)$ on $\mathcal{S}^\downarrow$ indexed by $1 \leq n_1 \leq n - 1$.*

PROOF. Let $\Pi$ be a homogeneous fragmentation with dislocation measure $\nu$. The result is obtained by conditioning on the size of the block $B_1$ containing 1. We (have to) take $g(n, n_1)$ as the total rate associated with the formation of a particular block $B_1$ of $n_1$ out of $n$ elements. Then $\binom{n-1}{n_1-1} g(n, n_1) = \Phi(n-1 : n - n_1)$ as in (7), so that

$$\Phi(n-1) = \sum_{n_1=1}^{n-1} \binom{n-1}{n_1-1} g(n, n_1) = \mathbb{P}_\nu(\Pi_n \neq \{[n]\})$$

$$= \int_{\mathcal{S}^\downarrow} \left(1 - \sum_{j=1}^\infty s_j^n\right) \nu(d\mathbf{s}), \tag{10}$$

as in (8), is the total rate of formation of partitions of $[n]$ with at least 2 parts. Then $p_{\widehat{\nu}(n,n_1)}(n_2, \ldots, n_k)$ is the conditional probability, given the particular set $B_1$, that the remaining $n - n_1$ elements are partitioned as they must be to make a particular partition of $[n]$ into blocks of sizes $n_1, \ldots, n_k$. To be more precise, we can take

$$\widehat{\nu}(n, n_1)(d\mathbf{s}) = \frac{1}{g(n, n_1)} \int_{\mathcal{S}^\downarrow} \sum_{i \geq 1} r_i^{n_1} (1 - r_i)^{n-n_1} \delta_{\hat{\mathbf{r}}_i/(1-r_i)}(d\mathbf{s}) \nu(d\mathbf{r}),$$

where $\hat{\mathbf{r}}_i$ is the vector $\mathbf{r}$ with component $r_i$ omitted. By Kingman's paintbox representation and conditioning on the color $i$ of the first block, we then get for all partitions with block sizes $(n_1, \ldots, n_k)$ in order of least element

$$p_{\widehat{\nu}(n,n_1)}(n_2, \ldots, n_k)$$
$$= \int_{\mathcal{S}^\downarrow} p_{\mathbf{s}}(n_2, \ldots, n_k) \widehat{\nu}(n, n_1)(d\mathbf{s})$$
$$= \frac{1}{g(n, n_1)} \int_{\mathcal{S}^\downarrow} \sum_{i \geq 1} r_i^{n_1} (1 - r_i)^{n-n_1} p_{\hat{\mathbf{r}}_i/(1-r_i)}(n_2, \ldots, n_k) \nu(d\mathbf{r})$$
$$= \frac{1}{g(n, n_1)} p_\nu(n_1, \ldots, n_k),$$

where by convention $p_\mathbf{s} = p_{\delta_\mathbf{s}}$. □

This discussion simplifies greatly for measures $\nu$ with the special symmetry property introduced in the following definition:

DEFINITION 2. Let $\nu$ be a measure on $\mathcal{S}^\downarrow$, and let $\widehat{\nu}$ be a probability measure on $\mathcal{S}^\downarrow$. Say that $\nu$ *has $\widehat{\nu}$ as its factor*, if $\widehat{\nu}(n,n_1)$ in (9) can be chosen identically equal to $\widehat{\nu}$ for every $1 \leq n_1 < n$, that is,

$$p_\nu(n_1,\ldots,n_k) = g(n,n_1) p_{\widehat{\nu}}(n_2,\ldots,n_k) \tag{11}$$

for every composition $n_1,\ldots,n_k$ of $n \geq 2$ into at least 2 parts, and some function $g(n,n_1)$.

Note that $\nu$ may be sigma-finite, but that $\widehat{\nu}$ is always assumed to be a probability measure. It is obvious that if $\nu$ has factor $\widehat{\nu}$, then $\widehat{\nu}$ is unique. A rich class of measures $\nu$ which admit a factor $\widehat{\nu}$ is the class of Poisson–Dirichlet measures considered in the next section. It is an open problem [32], Problem 3.7, even for probability measures, to describe all measures $\nu$ on $\mathcal{S}^\downarrow$ which admit a factor $\widehat{\nu}$. Note that all binary dislocation measures trivially admit a factor, as well as ordered Dirichlet$(a,\ldots,a)$ including the Dirac mass at $(1/m,\ldots,1/m)$. The latter are just the remaining members of the Ewens–Pitman two-parameter family.

Following the formalism of [31], Corollary 13, given two random elements $V$ and $V'$ of $\mathcal{S}^\downarrow$, and a probability distribution $\widehat{\nu}$ on $\mathcal{S}^\downarrow$, say that $V'$ *is a $\widehat{\nu}$-fragmentation of* $V$ if the joint distribution of $V$ and $V'$ is the same as if $V'$ is derived from $V$ by shattering each fragment of $V$ independently in proportions determined by $\widehat{\nu}$.

THEOREM 6. *Let $\nu$ be a dislocation measure on $\mathcal{S}^\downarrow$, let $(\mathcal{T},\mu)$ be some self-similar CRT derived from fragmentation according to $\nu$, and let $\widehat{\nu}$ be a probability distribution on $\mathcal{S}^\downarrow$. Then the following two conditions are equivalent:*

  (i) *the measure $\nu$ has $\widehat{\nu}$ as a factor;*
  (ii) *the fine spinal mass partition of $(\mathcal{T},\mu)$ is a $\widehat{\nu}$-fragmentation of the coarse spinal mass partition of $(\mathcal{T},\mu)$.*

PROOF. According to Pitman ([31], Lemma 35), the fine spinal partition is a $\widehat{\nu}$-fragmentation of the coarse spinal partition if and only if, for all $n \geq 1$, in passing from the coarse spinal partition of $[n]$ generated by $T_{[n]}$ to the fine one, within each block of the coarse partition the fine partition is distributed according to $\mathbb{P}_{\widehat{\nu}}$, independently between blocks of the coarse partition. So fix some integer $n$ and let $B_1, \ldots, B_k$ be the blocks of the coarse spinal partition of $T_{[n+1]}$ restricted to $\{2, \ldots, n+1\}$, with respective sizes $n_1, \ldots, n_k$. Due to the fragmentation property of the trees $T_{[n]}, n \geq 1$, the corresponding fine spinal partition of $T_{[n+1]}$ is obtained by splitting independently $B_1$ according to $\mathbb{P}_{\widehat{\nu}(n+1,n+1-n_1)}$, $B_2$ according to $\mathbb{P}_{\widehat{\nu}(n+1-n_1,n+1-n_1-n_2)}, \ldots, B_k$ according to $\mathbb{P}_{\widehat{\nu}(n_k+1,1)}$, where $\widehat{\nu}(n+1,n+1-n_1), \widehat{\nu}(n+1-n_1,n+1-n_1-n_2), \ldots, \widehat{\nu}(n_k+1,1)$ are probability measures satisfying (9). The fine spinal partition is therefore a $\widehat{\nu}$-fragmentation of the coarse spinal partition if and only if $\widehat{\nu}(n,n_1)$ can be chosen equal to $\widehat{\nu}$ for all $1 \leq n_1 < n$. $\square$

**3. Poisson–Dirichlet fragmentations.** We now turn to a particular family of fragmentation processes, namely the *Poisson–Dirichlet fragmentations*, characterized by dislocation measures of type $\mathrm{PD}^*(\alpha,\theta)$, $0 < \alpha < 1$, $\theta > -2\alpha$, as defined below by (19). This family generalizes the family of previously studied *stable fragmentations* [29, 30], constructed from the *stable trees* $(\mathcal{T}_\beta, \mu_\beta)$ with index $\beta$, $1 < \beta < 2$. These stable CRTs were introduced and studied by Duquesne and Le Gall [14, 15] to which we refer for a rigorous construction. Roughly, $\mathcal{T}_\beta$ arises as the limit in distribution as $n \to \infty$ of rescaled critical Galton–Watson trees $T_n$, conditioned to have $n$ vertices, with edge-lengths $n^{1/\beta-1}$, and an offspring distribution $(\eta_k, k \geq 0)$ such that $\eta_k \sim Ck^{-1-\beta}$ as $k \to \infty$. It is endowed with a (random) probability measure $\mu_\beta$ which is the limit as $n \to \infty$ of the empirical measure on the vertices of $T_n$. Miermont [29] shows that the partition-valued process constructed by random sampling of leaves $L_1, L_2, \ldots$ from $(\mathcal{T}_\beta, \mu_\beta)$ according to $\mu_\beta$ (as explained at the end of Section 1.2) is a self-similar fragmentation with index $1/\beta - 1$, and dislocation measure $\nu_\beta$ defined for all nonnegative measurable function $f$ on $\mathcal{S}^\downarrow$ by

$$(12) \qquad \int_{\mathcal{S}^\downarrow} f(\mathbf{s})\nu_\beta(d\mathbf{s}) = \frac{\beta^2 \Gamma(2-1/\beta)}{\Gamma(2-\beta)} \mathbb{E}\bigg[Tf\bigg(\frac{\Delta_1}{T}, \frac{\Delta_2}{T}, \ldots\bigg)\bigg]$$

(and no erosion). Here $T = \sum_{i=1}^{\infty} \Delta_i$ where $\Delta_1 > \Delta_2 > \cdots$ are the points of a Poisson process on $(0,\infty)$ with intensity $(\beta\Gamma(1-1/\beta))^{-1} x^{-1/\beta-1} \, dx$. Besides, cutting the stable tree $\mathcal{T}_\beta$ at nodes (see [30]), Miermont obtained a self-similar fragmentation with index $1/\beta$ and the same dislocation measure $\nu_\beta$.

3.1. *Definition and factorization property.* For $0 \leq \alpha < 1, \theta > -\alpha$, let $\mathrm{PD}(\alpha,\theta)$ denote the two-parameter Poisson–Dirichlet distribution on $\mathcal{S}^\downarrow$, defined as the distribution of the decreasing rearrangement of its size-biased presentation, which is

$$(13) \qquad W_1, (1-W_1)W_2, (1-W_1)(1-W_2)W_3, \ldots$$

for $W_i$, which are independent beta$(1-\alpha, i\alpha+\theta)$ variables. The formula for the corresponding EPPF is [32], Theorem 3.2,

$$(14) \quad p_{\mathrm{PD}(\alpha,\theta)}(n_1, \ldots, n_k) = \frac{\alpha^{k-1}[1+\theta/\alpha]_{k-1}}{[1+\theta]_{n-1}} \prod_{i=1}^{k} [1-\alpha]_{n_i-1}$$

for every composition $(n_1, \ldots, n_k)$ of $n$, where $[x]_n = \Gamma(x+n)/\Gamma(x)$ is a rising factorial. It is evident by inspection of this formula and (11) that the probability measure $\mathrm{PD}(\alpha,\theta)$ admits $\mathrm{PD}(\alpha, \theta+\alpha)$ as a factor for every $0 < \alpha < 1$ and $\theta > -\alpha$. Following Miermont [29] we now consider the rescaled measure

$$(15) \quad \mathrm{PD}^*(\alpha,\theta) := \frac{\Gamma(1+\theta/\alpha)}{\Gamma(1+\theta)} \mathrm{PD}(\alpha,\theta),$$

which is defined in the first instance for $0 < \alpha < 1$ and $-\alpha < \theta$. It is known ([32], Corollary 3.9) that for $0 < \alpha < 1$ there is the absolute continuity relation

$$(16) \quad \mathrm{PD}^*(\alpha,\theta)(d\mathbf{s}) = (S_\alpha(\mathbf{s}))^{\theta/\alpha} \mathrm{PD}(\alpha,0)(d\mathbf{s}),$$

where $S_\alpha(\mathbf{s})$ is the $\alpha$-*diversity* which is almost surely associated to a sequence $\mathbf{s} = (s_1, s_2, \ldots)$ with distribution $\mathrm{PD}(\alpha,0)$ by the formula

$$(17) \quad S_\alpha(\mathbf{s}) := \Gamma(1-\alpha) \lim_{j \to \infty} j s_j^\alpha.$$

The $\mathrm{PD}(\alpha,\theta)$ distribution is recovered from (16) for $-\alpha < \theta$ by normalization as in (15). The $\alpha$-diversity $S_\alpha$, which has a Mittag–Leffler distribution (see, e.g., [32], (0.43)), appears variously disguised in different contexts, for example, as a local time variable ([32], page 10), or again as $S_\alpha = T^{-\alpha}$ for a positive stable variable $T$ of index $\alpha$. Indeed, if such a $T$ is constructed as $T = \sum_{i=1}^{\infty} \Delta_i$ where $\Delta_1 > \Delta_2 > \cdots$ are the points of a Poisson process on $(0,\infty)$ with intensity $\alpha(\Gamma(1-\alpha))^{-1} x^{-\alpha-1} dx$, then

$$(\Delta_1/T, \Delta_2/T, \ldots) =_d \mathrm{PD}(\alpha, 0)$$

and, according to [32], (4.45),

$$S_\alpha(\Delta_1/T, \Delta_2/T, \ldots) = T^{-\alpha} \quad \text{a.s.},$$

so that for every nonnegative measurable function $f$ of $\mathbf{s} = (s_1, s_2, \ldots) \in \mathcal{S}^\downarrow$,

$$(18) \quad \int_{\mathcal{S}^\downarrow} f(\mathbf{s}) \mathrm{PD}^*(\alpha,\theta)(d\mathbf{s}) = \mathbb{E}[T^{-\theta} f(\Delta_1/T, \Delta_2/T, \ldots)].$$

LEMMA 7. *For each $0 < \alpha < 1$, let $\mathrm{PD}^*(\alpha,\theta)$ be the measure defined on $\mathcal{S}^\downarrow$ for each real $\theta$ by either (16) or (18). Then for $-2\alpha < \theta$, this measure* $\mathrm{PD}^*(\alpha, \theta)$ *is also the unique measure with no mass at* $(1, 0, 0, \ldots)$ *whose EPRF is given for* $k \geq 2$ *by*

$$(19) \qquad p_{\mathrm{PD}^*(\alpha,\theta)}(n_1, \ldots, n_k) = \frac{\alpha^{k-1}\Gamma(k+\theta/\alpha)}{\Gamma(n+\theta)} \prod_{i=1}^{k} [1-\alpha]_{n_i-1}$$

*and for* $k = 1$ *by the same formula for* $-\alpha < \theta$, *and by* $\infty$ *for* $-2\alpha < \theta \leq -\alpha$. *Basic integrability properties of this extended family of Poisson–Dirichlet measures are*

$$(20) \qquad \int_{\mathcal{S}^\downarrow} \mathrm{PD}^*(\alpha,\theta)(d\mathbf{s}) < \infty \quad \Leftrightarrow \quad \theta > -\alpha;$$

$$(21) \qquad \int_{\mathcal{S}^\downarrow} (1-s_1) \mathrm{PD}^*(\alpha,\theta)(d\mathbf{s}) < \infty \quad \Leftrightarrow \quad \theta > -2\alpha.$$

*For each choice of* $(\alpha, \theta)$ *with* $\theta > -2\alpha$ *the measure* $\mathrm{PD}^*(\alpha, \theta)$ *has the probability distribution* $\mathrm{PD}(\alpha, \theta + \alpha)$ *as its factor.*

PROOF. Following Miermont ([29], Section 3.3) we observe from (14) and (15) that the formula (19) holds in the first instance for all $\theta > -\alpha$, and that the right-hand side of (19) is analytic in $\theta$ for $\mathrm{Re}(\theta) > -2\alpha$, when $k \geq 2$. To get (19) for all $\theta > -2\alpha$, note that the left-hand side of (19) can be written as $\mathbb{E}[T^{-\theta}Y]$ where $Y$ is some positive r.v. depending on $n_1, \ldots, n_k$ and then

$$\mathbb{E}[T^{-\theta}Y] = \mathbb{E}[T^{-\theta}Y \mathbf{1}_{\{T<1\}}] + \mathbb{E}[T^{-\theta}Y \mathbf{1}_{\{T \geq 1\}}],$$

where the first term is finite for all $\theta \in \mathbb{R}$, hence an entire function of $\theta$. So the second term for $\theta > -\alpha$ equals a function that is analytic for $\mathrm{Re}(\theta) > -2\alpha$. We claim that this equality of functions extends to $\theta > -2\alpha$. Indeed, consider some nonnegative r.v. $Z$ such that $M(t) := \mathbb{E}[e^{tZ}] < \infty$ for $t < r_1$ and $M(t) = N(t)$ for $t < r_1$ where $N$ is analytic for $\mathrm{Re}(t) < r_2$ with $0 < r_1 < r_2$. Then the identity for $t < r_1$ gives the power series expansion $N(t) = \sum_0^\infty t^n \mathbb{E}[Z^n]/n!$ for $|t| < r_1$. Since $N$ is analytic for $\mathrm{Re}(t) < r_2$, this power series converges and this identity holds also for $|t| < r_2$. Hence for $0 \leq t < r_2$ we can compute by monotone convergence $M(t) = \sum_0^\infty t^n \mathbb{E}[Z^n]/n! = N(t)$. Hence (19).

The fact (20) comes from formula (0.40) and the following line in [32]. As for (21), we have seen in Section 2.2 that this integrability condition holds if and only if the expressions in (19) are finite for every choice of $n_1, \ldots, n_k$ with $k \geq 2$, and this is clear by inspection of (19). □

The infinite measure $\mathrm{PD}^*(\alpha, -\alpha)$ was already used and studied by Basdevant [7] in the context of Ruelle's probability cascades.

REMARKS. I. For $0 < \alpha < 1, \theta > -\alpha$, the EPPF (14) gives

$$\mathbb{P}_{\mathrm{PD}(\alpha,\theta)}(\Pi_n \neq \{[n]\}) = 1 - \frac{[1-\alpha]_{n-1}}{[1+\theta]_{n-1}} \tag{22}$$

and hence

$$\mathbb{P}_{\mathrm{PD}^*(\alpha,\theta)}(\Pi_n \neq \{[n]\}) = \frac{\Gamma(1+\theta/\alpha)}{\Gamma(1+\theta)}\left(1 - \frac{[1-\alpha]_{n-1}}{[1+\theta]_{n-1}}\right) \tag{23}$$

in the first instance for $0 < \alpha < 1, \theta > -\alpha$, and then by analytic continuation for $0 < \alpha < 1, \theta > -2\alpha$, with values of the right-hand side defined by continuity for $\theta = -\alpha$ or $\theta = -1$. To see that the left-hand side of (23) is analytic in this range, observe that for each $n$ this function of $(\alpha, \theta)$ is just a finite sum of the functions in (19) weighted by combinatorial coefficients.

II. From the fact (13) that a size-biased pick from $\mathrm{PD}(\alpha, \theta)$ has beta$(1-\alpha, \alpha+\theta)$ distribution for $0 < \alpha < 1, \theta > -\alpha$, we can write down

$$s\sum_{j=1}^{\infty} \mathrm{PD}(\alpha,\theta)(s_j \in ds) = \frac{\Gamma(1+\theta)}{\Gamma(1-\alpha)\Gamma(\alpha+\theta)} s^{-\alpha}(1-s)^{\alpha+\theta-1}\, ds$$

$$(0 < s < 1)$$

and hence for $-2\alpha < \theta$ by analytic continuation

$$s\sum_{j=1}^{\infty} \mathrm{PD}^*(\alpha,\theta)(s_j \in ds) = \frac{\alpha\Gamma(2+\theta/\alpha)}{\Gamma(1-\alpha)\Gamma(1+\alpha+\theta)} s^{-\alpha}(1-s)^{\alpha+\theta-1}\, ds \tag{24}$$

$$(0 < s < 1).$$

The image of this measure by the change of variable $x = -\log s$ is the corresponding Lévy measure

$$\Lambda_{\alpha,\theta}(dx) = \frac{\alpha\Gamma(2+\theta/\alpha)}{\Gamma(1-\alpha)\Gamma(1+\alpha+\theta)} e^{-x(1-\alpha)}(1-e^{-x})^{\alpha+\theta-1}\, dx \tag{25}$$

$$(0 < x < \infty).$$

From Theorem 6 we now deduce:

COROLLARY 8. *For each $0 < \alpha < 1, \theta > -2\alpha$, let $(\mathcal{T}_{\alpha,\theta}, \mu)$ be some CRT derived from fragmentation process with dislocation measure $\mathrm{PD}^*(\alpha, \theta)$. The sequence of discrete fragmentation trees $(T_{[n]}, n \geq 1)$ embedded in $(\mathcal{T}_{\alpha,\theta}, \mu)$ is governed by fragmentations of $[n]$ according the EPPF obtained by normalization of formula (19) by formula (23). The fine spinal mass partition of $(\mathcal{T}_{\alpha,\theta}, \mu)$ is a $\mathrm{PD}(\alpha, \alpha+\theta)$-fragmentation of the coarse spinal mass partition* of $(\mathcal{T}_{\alpha,\theta},\mu)$, which is derived from the range of $1 - e^{-\xi}$ for the pure jump subordinator $\xi$ with Lévy measure (25) and Laplace exponent

$$(26)\quad \Phi_{\alpha,\theta}(z) = \begin{cases} \dfrac{\alpha\Gamma(2+\theta/\alpha)}{(\alpha+\theta)\Gamma(1-\alpha)}\bigg(\dfrac{(1+\theta)\Gamma(1-\alpha)}{\Gamma(2+\theta)} \\ \qquad\qquad - \dfrac{(z+1+\theta)\Gamma(z+1-\alpha)}{\Gamma(z+2+\theta)}\bigg), & \theta \neq -\alpha, \\ \dfrac{\alpha}{\Gamma(1-\alpha)}\bigg(\dfrac{\Gamma'(z+1-\alpha)}{\Gamma(z+1-\alpha)} - \dfrac{\Gamma'(1-\alpha)}{\Gamma(1-\alpha)}\bigg), & \theta = -\alpha. \end{cases}$$

Last, for $\theta \in (-2\alpha, -\alpha)$ we have an interesting regime where Proposition 4 applies along with the asymptotic theory of consistent Markov branching models in [26]. Specifically,

COROLLARY 9. *For $\theta \in (-2\alpha, -\alpha)$, let $(T_{[n]}, n \geq 1)$ be a Markov branching model derived from a self-similar fragmentation with dislocation measure $\mathrm{PD}^*(\alpha,\theta)$. Adding unit edge lengths to $T_{[n]}$, there is the convergence in probability*

$$(27)\quad \frac{|\alpha+\theta|\Gamma(1-\alpha)}{\alpha\Gamma(2+\theta/\alpha)} \times \frac{T_{[n]}}{n^{|\theta+\alpha|}} \to \mathcal{T}_{(\theta+\alpha,\mathrm{PD}^*(\alpha,\theta))}$$

*for the Gromov–Hausdorff topology, where the limit is a self-similar fragmentation CRT of index $\theta+\alpha$ and dislocation measure $\mathrm{PD}^*(\alpha,\theta)$.*

PROOF. Note from (24) that

$$\mathrm{PD}^*(\alpha,\theta)(s_1 \leq 1-\varepsilon) \sim \frac{\alpha\Gamma(2+\theta/\alpha)}{|\alpha+\theta|\Gamma(1-\alpha)\Gamma(1+\alpha+\theta)}\varepsilon^{\alpha+\theta} \qquad \text{as } \varepsilon \downarrow 0.$$

Then Theorem 2 of [26] applies [$\Lambda_{\alpha,\theta}$ clearly also satisfies $\int^\infty x^\rho \Lambda_{\alpha,\theta}(d\mathbf{s}) < \infty$ for some $\rho > 0$], which gives (27). □

3.2. *Stable fragmentations.* The case $1/2 < \alpha < 1$ is of special interest. Then $-2\alpha < -1 < -\alpha$, so we can take $\theta = -1$ in (24), and then the Lévy measure (25) is of the form

$$(28)\quad \Lambda(dx) = c_b(1-e^{-x})^{-b-1}e^{-bx}\,dx$$

for some constant $c_b > 0$ and $b = 1 - \alpha$. It is known [19] that if $\xi$ is a subordinator with this Lévy measure, for any $b \in (0,1)$, then the closure of the range of $e^{-\xi}$ is reversible and identical in law with the zero set of a Bessel bridge of dimension $2 - 2b$. The corresponding distribution of ranked lengths of intervals is then known to be $\mathrm{PD}(b,b)$ ([32], Corollary 4.9). Miermont ([29], page 444) found the same Lévy measure, up to a scaling constant, for the subordinator associated with the self-similar fragmentation of index $\alpha - 1 \in (-1/2, 0)$ that he derived from the stable CRT $\mathcal{T}_\beta$ of index $\beta = 1/\alpha \in (1, 2)$. Here we have reversed this line of reasoning, and constructed $\mathcal{T}_\beta$ directly from combinatorial considerations, without relying on the relation between the height process of $\mathcal{T}_\beta$ and the stable process of index $\beta$, which was the basis of the work of Duquesne and Le Gall [14, 15]. As byproducts of this argument, we have a number of refinements of earlier work on $\mathcal{T}_\beta$, which we summarize in the following corollary of previous results.

COROLLARY 10. *For each $\alpha \in (1/2, 1)$, corresponding to $\beta = 1/\alpha \in (1, 2)$ the dislocation measure $\mathrm{PD}^*(\alpha, -1)$ derived from the two-parameter Poisson–Dirichlet family as in (19) has $\mathrm{PD}(\alpha, \alpha - 1)$ as a factor. Let $(T_{[n]}, n = 1, 2 \ldots)$ be a consistent family of combinatorial trees governed by fragmentation according to $\mathrm{PD}^*(\alpha, -1)$. Then:*

1. *The tree $T_{[n]}$ is identical in law to the combinatorial tree with $n$ leaves derived by sampling according the mass measure in the stable tree $\mathcal{T}_\beta$ of index $\beta$, and $\mathcal{T}_\beta$ may be constructed from the sequence of combinatorial trees $(T_{[n]}, n \geq 1)$, as indicated in [26], Theorem 2, or Corollary 3.*
2. *The distribution of the coarse spinal mass partition of $\mathcal{T}_\beta$ is $\mathrm{PD}(1 - \alpha, 1 - \alpha)$.*
3. *The coarse spinal interval partition of $[0, 1]$ derived from $\mathcal{T}_\beta$ is exchangeable, with the same distribution as the collection of excursion intervals of a Bessel bridge of dimension $2\alpha$. The $(1 - \alpha)$-diversity of this interval partition—defined in a way similar to (17)—is a multiple of the height of a leaf picked at random from the mass measure of $\mathcal{T}_\beta$. This height has the same tilted Mittag–Leffler distribution as a multiple of the local time at $0$ of the Bessel bridge of dimension $2\alpha$.*
4. *The corresponding fine spinal mass partition of $\mathcal{T}_\beta$ is a $\mathrm{PD}(\alpha, \alpha - 1)$-fragmentation of the coarse spinal mass partition.*
5. *The unconditional distribution of the fine spinal mass partition of $\mathcal{T}_\beta$ is $\mathrm{PD}(\alpha, 1 - \alpha)$.*
6. *The conditional distribution of the coarse spinal mass partition of $\mathcal{T}_\beta$ given the fine one is provided by the operator of $\mathrm{PD}(\gamma, \gamma)$ coagulation, as defined in [31], for $\gamma = (1 - \alpha)/\alpha$.*
7. *Conditionally given the fine spinal mass partition of $\mathcal{T}_\beta$, the corresponding collection of subtrees obtained by removing the spine, modulo isomorphism and rescaling trees $T$ of mass $m$ to $m^{-(1-\alpha)} T$, is a collection of independent copies of $\mathcal{T}_\beta$.*

PROOF. All but items 3, 5 and 6 follow immediately from the previous development. Items 5 and 6 are read from items 2 and 4 by the more general coagulation/fragmentation duality relation for the PD family provided by [31], Theorem 12.

The first assertion of item 3 is obvious from previous discussion. To get the remaining assertions of item 3, denote by $\Pi$ the coarse spinal partition of $\mathbb{N} \setminus \{1\}$ derived from $\mathcal{T}_\beta$ and for each $n$ by $K_n$ the number of blocks of $\Pi_n$. On the one hand, Theorem 3.8 and Lemma 3.13 of [32] ensure that $n^{\alpha-1} K_n$ converges a.s. to the $(1-\alpha)$-diversity, that moreover has a tilted Mittag–Leffler distribution, given precisely by formula (3.27) of [32]. On the other hand, $K_n + 1$ is the height of leaf 1 in $T_{[n]}$. Hence by Corollary 9, $n^{\alpha-1} K_n$ converges in probability to a multiple of the height of a leaf taken at random from the mass measure of $\mathcal{T}_\beta$. This gives the second assertion of item 3. For the last assertion of item 3, recall from Theorem 5.3 in [33], that the decreasing sequence of excursion lengths of a Bessel bridge of dimension $2\alpha$ is absolutely continuous with respect to the distribution of the decreasing sequence of excursion lengths until time 1 of a Bessel process of dimension $2\alpha$. By [32], Section 4.4, the local time at level 0 until time 1 of this Bessel process is proportional to the $(1-\alpha)$-diversity constructed from the sequence of lengths of its excursions. Hence a similar result holds for the Bessel bridge. $\square$

For more information about the distribution of random partitions in the PD family, see [20] and [34]. In the limiting case when $\beta \uparrow 2$, the above results reduce to the description of the interval partition derived from the spinal decomposition of the Brownian CRT, which is well known to be distributed like the partition generated by excursions of a Brownian bridge. See [5] for applications of this decomposition to the asymptotics of random mappings. The structure of the fine spinal partition of $\mathcal{T}_\beta$ for $1 < \beta < 2$ has no analogue for $\beta = 2$, because in the Brownian tree all splits are binary.

**4. Invariance under uniform re-rooting.** It is of particular interest to consider fragmentation trees with additional symmetry properties. A well-known property of the stable tree $\mathcal{T}_\beta$ with index $\beta \in (1, 2]$, established by Aldous [2] for the Brownian CRT with $\beta = 2$, and by Duquesne and Le Gall [15], Proposition 4.8, for $\beta \in (1, 2)$, is *invariance under uniform re-rooting*. See also [13]. Let us first introduce the discrete analogue of this property.

For a tree $T_{[n]}$ with leaves labeled by $[n]$, let $T_{[n]}^{(\text{ROOT} \leftrightarrow 1)}$ denote the tree with leaves labeled by $[n]$ obtained by re-rooting $T_{[n]}$ at 1 and relabeling the original root by 1. See, for instance, Figure 3.

DEFINITION 3. Let $(T_{[n]}, n \geq 1)$ be a consistent Markov branching model. We say that the Markov branching model is *invariant under uniform re-rooting* if for all $n \geq 1$

$$T_{[n]} \stackrel{\text{law}}{=} T_{[n]}^{(\text{ROOT} \leftrightarrow 1)}.$$

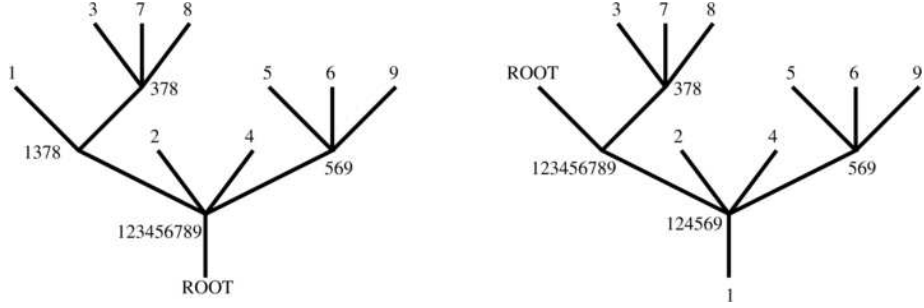

FIG. 3. *A fragmentation tree $T_{[9]}$ and its re-rooted counterpart $T_{[9]}^{(\mathrm{ROOT}\leftrightarrow 1)}$.*

Note that due to the exchangeability of leaf labels, leaf 1 is indeed a uniformly picked leaf of the de-labeled combinatorial tree shape. Due to the exchangeability of leaf labels, invariance under uniform re-rooting is in fact a property of de-labeled combinatorial tree shapes.

DEFINITION 4. Let $(\mathcal{T},\mu)$ be a CRT rooted at $\rho$ and conditionally on $(\mathcal{T},\mu)$, let $(L_1, L_2, \ldots)$ be a sample of leaves i.i.d. with distribution $\mu$. Let then $\mathcal{T}^{[L_1]}$ denote the tree $\mathcal{T}$ re-rooted at $L_1$. We say that $(\mathcal{T},\mu)$ is invariant under uniform re-rooting if for all $n \geq 1$, the law of the reduced subtree $\mathcal{R}(\mathcal{T}, L_1, \ldots, L_n)$ of $\mathcal{T}$ spanned by the root $\rho$ and $L_1, \ldots, L_n$ is invariant under re-rooting at $L_1$, that is,

$$\mathcal{R}(\mathcal{T}^{[L_1]}, \rho, L_2, \ldots, L_n) \stackrel{\mathrm{law}}{=} \mathcal{R}(\mathcal{T}, L_1, L_2, \ldots, L_n)$$

as an identity in law of combinatorial tree shapes with assignment of edge lengths.

Clearly, the invariance under uniform re-rooting of $(\mathcal{T},\mu)$ implies the invariance under uniform re-rooting of the sequence $(T_{[n]}, n \geq 1)$ of combinatorial trees associated with $(\mathcal{T},\mu)$. We will see that the converse is false [see the arguments after (36)].

REMARK. In [2, 13] a different formalism is used for the definition of invariance under uniform re-rooting, via height functions of ordered CRTs. Briefly, assuming that the CRT $(\mathcal{T},\mu)$ can be encoded into a continuous real-valued function $H$ on $[0,1]$, with $H(0) = H(1) = 0$, such that:

- $\mathcal{T}$ is isometric to the quotient space $([0,1], d_H)/\sim_H$, where
$$d_H(x,y) := H(x) + H(y) - 2 \min_{z \in [x,y]} H(z)$$
and
$$x \sim_H y \quad \Leftrightarrow \quad d_H(x,y) = 0$$

with the convention $[x, y] = [y, x]$, when $y < x$,
- $\mu$ is the measure induced by the projection of the Lebesgue measure on this quotient space,

then the invariance under uniform re-rooting is defined via $H^{[U]} \stackrel{\text{law}}{=} H$ where $U$ if uniformly distributed on $[0, 1]$ independently of $H$ and

$$(29) \quad H^{[u]}(x) := H(u) + H(u+x) - 2 \min_{z \in [u, x+u]} H(z), \qquad u, x \in [0, 1]$$

with the convention $u + x = u + x - 1$, when $u + x > 1$. It was proved in [25] that the structures of the combinatorial subtrees $\mathcal{R}(\mathcal{T}, L_1, \ldots, L_n)$, $n \geq 1$, derived from some self-similar fragmentation CRT $(\mathcal{T}, \mu)$ can be enriched with a consistent "uniform" order so as to encode the fragmentation CRTs into a continuous height function as described above, provided the dislocation measure is infinite. In that context, it is not hard to check that the height function definition and Definition 4 above are equivalent. Details are left to the reader.

The goal of this section is twofold: first to give a combinatorial proof, different from that given in [2, 13, 15], of the fact that the stable trees are invariant under uniform re-rooting; second to prove that among the self-similar fragmentation CRTs, the stable trees are the only ones, up to a scaling factor, to satisfy this invariance property.

For the Brownian CRT $(\mathcal{T}_2, \mu_2)$, we recall that the partition-valued process constructed by random sampling of leaves $L_1, L_2, \ldots$ according to $\mu_2$ is a self-similar fragmentation with index $a = -1/2$ and dislocation measure $\nu_2$ defined by $\nu_2(s_1 + s_2 \neq 1) = 0$ and

$$\nu_2(s_1 \in dx) = (2/\pi)^{1/2} x^{-3/2} (1-x)^{-3/2} \, dx, \qquad 1/2 \leq x < 1$$

(see [9]). The dislocation measure $\nu_\beta$ associated to the stable tree $\mathcal{T}_\beta$, $1 < \beta < 2$, is given by (12) and its self-similar index is $1/\beta - 1$.

THEOREM 11. (i) *[2, 15]. For all $\beta \in (1, 2]$, the stable tree $(\mathcal{T}_\beta, \mu_\beta)$ is invariant under uniform re-rooting.*

(ii) *Let $(\mathcal{T}, \mu)$ be a self-similar fragmentation CRT with parameters $(a, \nu)$ and suppose it is invariant under uniform re-rooting. Then there exists some $\beta \in (1, 2]$ and some constant $C > 0$ such that $\nu = C \nu_\beta$ and $a = 1/\beta - 1$.*

REMARK. According to [13], a stronger invariance result is available for the height functions $H$ of stable trees (and more generally Lévy trees), which is that $H^{[u]}$, as defined in (29), is distributed as $H$ for each fixed $u \in [0, 1]$. See also [28] for the Brownian CRT.

The rest of this section is devoted to the proof of Theorem 11.

4.1. *Spinal decomposition and proof of Theorem 11*(i). The first step is to consider the spinal decomposition of trees invariant under uniform re-rooting: one consequence of this invariance is that the coarse spinal interval partition of $[0,1]$ derived from the tree is reversible (in fact an exchangeable interval partition of $[0,1]$; see [19]). The class of trees with this property is significantly restricted by the following proposition.

PROPOSITION 12. *Let $(T_{[n]}, n \geq 1)$ be a sequence of combinatorial trees associated with some self-similar fragmentation CRT $(\mathcal{T}, \mu)$ with dislocation measure $\nu$, and let $\xi$ be the subordinator describing the evolution of the mass fragment containing 1 in an associated homogeneous fragmentation process (cf. Section 2.1).*

*(i) The coarse spinal composition $\mathcal{C}_n$ of $n$ derived from $T_{[n+1]}$ [as defined in (3)] is reversible for each $n$ if and only if $\xi$ has a Lévy measure of the form*

$$\Lambda(dx) = c(1 - e^{-x})^{-b-1} e^{-bx} \, dx \tag{30}$$

*for some $0 < b < 1$ and some constant $c > 0$.*

*(ii) There cannot exist a self-similar fragmentation CRT with a Lévy measure of this form when $b > 1/2$.*

PROOF. Part (i) is read from [19], Theorem 10.1, just using that $(\mathcal{C}_n, n \geq 1)$ is a regenerative composition structure. For part (ii), from (5)

$$\Lambda(dx) = e^{-x} \sum_{i \geq 1} \nu(-\log s_i \in dx), \qquad x > 0,$$

and (30) we deduce by the transformation $z = e^{-x}$, $x = -\log(z)$, $dx = -dz/z$

$$\sum_{i \geq 1} \nu(s_i \in dz) = c(1-z)^{-b-1} z^{b-2} \, dz, \qquad z \in (0,1).$$

Since $\nu$ is supported by decreasing sequences with $\sum_{i=1}^{\infty} s_i = 1$, $s_i \leq 1/i$ for all $i \geq 1$. In particular,

$$\nu(s_1 \in dz) = c(1-z)^{-b-1} z^{b-2} \, dz, \qquad z \in (1/2, 1). \tag{31}$$

Using the fact that for $z \in (0, 1/2)$

$$z^{-b}(1-z)^{b-2} > (1-z)^{-b-1} z^{b-1}$$
$$\iff (1-z)^{2b-1} > z^{2b-1}$$
$$\iff b > 1/2,$$

we see that for $b > 1/2$

$$\int_{(0,1)} (1-z)\nu(s_1 \in dz) \geq c \int_{(1/2,1)} (1-z)^{-b} z^{b-2} dz$$
$$= c \int_{(0,1/2)} z^{-b} (1-z)^{b-2} dz$$
$$> c \int_{(0,1/2)} (1-z)^{-b-1} z^{b-1} dz$$
$$\geq \sum_{i \geq 2} \int_{(0,1)} z\nu(s_i \in dz)$$

by (31). On the other hand, we have

$$\int_{(0,1)} (1-z)\nu(s_1 \in dz) = \int_{(0,1)} z\nu\left(\sum_{i \geq 2} s_i \in dz\right)$$
$$= \int_{\mathcal{S}^{\downarrow}} \sum_{i \geq 2} s_i \nu(d\mathbf{s})$$
$$= \sum_{i \geq 2} \int_{\mathcal{S}^{\downarrow}} s_i \nu(d\mathbf{s})$$
$$= \sum_{i \geq 2} \int_{(0,1)} z\nu(s_i \in dz),$$

which contradicts the inequality obtained in the preceding calculation. $\square$

The Lévy measure associated with some fragmentation tree invariant under uniform re-rooting is therefore of the form (30) for some $0 < b \leq 1/2$. We recall that the Lévy measures associated to $\beta$-stable trees are of this form for $b = 1 - 1/\beta$ (see Section 3.2 for $1 < \beta < 2$ and [9] for $\beta = 2$) which covers the range $(0, 1/2]$ when $\beta$ varies in $(1, 2]$.

PROOF OF THEOREM 11(i). Let $(\mathcal{T}_\beta, \mu_\beta)$ be some stable CRT with index $\beta \in (1, 2]$. According to the previous proposition, its coarse spinal interval partition of [0,1] is reversible. Items 4 and 7 of Corollary 10 then give a construction of this CRT from its coarse spinal interval partition, via its fine spinal mass partition, that ensures the invariance under re-rooting property. $\square$

4.2. *Characterization of the dislocation measure and proof of Theorem 11(ii).* In general the Lévy measure does not characterize the dislocation measure of the fragmentation tree, that is, two different dislocation measures may lead to the same Lévy measure $\Lambda$; see Haas [23] for an example.

However, this complication no longer arises when the set of fragmentation trees is restricted to the ones invariant under uniform re-rooting.

PROPOSITION 13. *Let $(\mathcal{T}, \mu)$ be a self-similar fragmentation CRT with parameters $(a, \nu)$ and suppose it is invariant under uniform re-rooting. Then the dislocation measure $\nu$ can be re-constructed from the Lévy measure $\Lambda$ associated to the tagged fragment.*

Together with Proposition 12, this implies that:

COROLLARY 14. *The dislocation measure of a self-similar fragmentation CRT invariant under uniform re-rooting is proportional to $\nu_\beta$ for some $\beta \in (1, 2]$.*

In order to prove Proposition 13, we first set up two lemmas. In the rest of this section, the CRT $(\mathcal{T}, \mu)$ with parameters $(a, \nu)$ is fixed and supposed to be invariant under uniform re-rooting. A sample of leaves $L_i, i \geq 1$, is given and we consider the associated partition-valued fragmentation $\Pi$. We call $p_n$ the probabilities

$$p_n(n_1, \ldots, n_k) = \mathbb{P}(\Pi_n(t_n) = \{\{1, \ldots, n_1\}, \{n_1 + 1, \ldots, n_1 + n_2\}, \ldots,$$
$$\{n_1 + \cdots + n_{k-1} + 1, \ldots, n\}\})$$
$$= \frac{p_\nu(n_1, \ldots, n_k)}{\Phi(n-1)},$$

where $t_n$ is the first time when $\Pi_n$ differs from $[n]$ and $(n_1, \ldots, n_k)$ denotes any composition of $n$ with $k \geq 2$ (in other words, the probabilities $p_n$ are the EPPFs obtained by conditioning $\mathbb{P}_\nu$ on $\{\Pi_n \neq \{[n]\}\}$ in the proof of Lemma 5). Note in particular that

$$(32) \qquad \sum_{k=2}^{n} \sum_{(n_1, \ldots, n_k)} \frac{n!}{n_1! \cdots n_k!} \frac{1}{k!} p_n(n_1, \ldots, n_k) = 1,$$

where the sum is over all compositions of $n$; see [32], Exercise 2.1.3.

LEMMA 15. *For all compositions $(n_1, \ldots, n_k)$ of $n$ with $k \geq 2$*

$$(33) \quad \begin{aligned} & p_n(n_1, \ldots, n_k) p_{n_1}(1, n_1 - 1) \\ & = p_n(n_2 + \cdots + n_k + 1, n_1 - 1) p_{n-n_1+1}(1, n_2, \ldots, n_k) \end{aligned}$$

*with the convention, when $n_1 = 1$, that the probabilities involving expressions with a term $n_1 - 1 = 0$ are all equal to 1.*

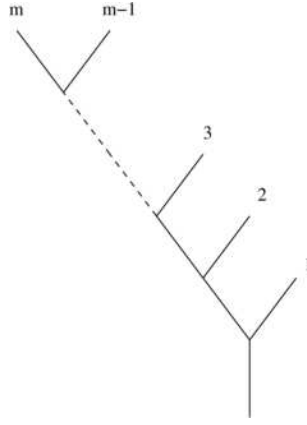

FIG. 4. *This configuration always happen with positive probability.*

PROOF. Consider the following fragmentation scheme: the first time at which the block $\{1,\ldots,n\}$ splits, it splits in blocks $\{1,\ldots,n_1\}$, $\{n_1+1,\ldots,n_1+n_2\},\ldots,\{n_1+\cdots+n_{k-1}+1,\ldots,n\}$; then the first of these blocks splits in $\{1\}$, $\{2,\ldots,n_1\}$. We are not really interested in the further evolution of $\{2,\ldots,n_1\}$, $\{n_1+1,\ldots,n_1+n_2\},\ldots,\{n_1+\cdots+n_{k-1}+1,\ldots,n\}$; let us just say that it is in a configuration which happens with a (strictly) positive probability, say $r_n(n_1,\ldots,n_k)$ (e.g., evolutions as in Figure 4). Consider then the discrete tree with leaf labels obtained from this fragmentation scheme. The probability that the tree with $n$ leaves $\mathcal{R}(\mathcal{T},L_1,L_2,\ldots,L_n)$ has this labeled shape is exactly

(34) $$p_n(n_1,\ldots,n_k)p_{n_1}(1,n_1-1)r_n(n_1,\ldots,n_k).$$

Now, look at the same tree rooted at $L_1$, that is, $\mathcal{R}(\mathcal{T}^{L_1},\rho,L_2,\ldots,L_n)$; cf. Figure 5. Starting from the root $L_1$, the corresponding fragmentation scheme evolves as follows: $\{\rho,2,\ldots,n\}$ first splits in $\{2,\ldots,n_1\}$, $\{\rho,n_1+1,\ldots,n\}$. Then $\{\rho,n_1+1,\ldots,n\}$ splits in $\{\rho\},\{n_1+1,\ldots,n_1+n_2\},\ldots,\{n_1+\cdots+n_{k-1}+1,\ldots,n\}$. And the blocks $\{2,\ldots,n_1\}$, $\{n_1+1,\ldots,n_1+n_2\},\ldots,\{n_1+\cdots+n_{k-1}+1,\ldots,n\}$ then all split according to the same configuration as in the previous scheme. By invariance under uniform re-rooting, the subtree $\mathcal{R}(\mathcal{T}^{[L_1]},\rho,L_2,\ldots,L_n)$ is distributed as $\mathcal{R}(\mathcal{T},L_1,L_2,\ldots,L_n)$, and therefore, the probability that $\mathcal{R}(\mathcal{T}^{[L_1]},\rho,L_2,\ldots,L_n)$ has this labeled shape is

(35) $p_n(n_1-1,n_2+\cdots+n_k+1)p_{n-n_1+1}(1,n_2,\ldots,n_k)r_n(n_1,\ldots,n_k).$

By invariance under re-rooting, the probabilities in (34) and (35) are equal. This yields (33), since $r_n(n_1,\ldots,n_k)\neq 0$. $\square$

REMARK. It is easy to check that the probabilities $p_n$ associated to the stable trees, which are obtained by normalization of formula (19) by (23) with $\theta = -1$, satisfy relations (33).

LEMMA 16. *The probabilities $p_3(2,1)$, $p_3(1,1,1)$ and $p_n(1, n-1)$, $\forall n \geq 2$, are determined by the Lévy measure $\Lambda$.*

PROOF. Consider $\Pi^0$, the homogeneous fragmentation constructed from $\Pi$ by time-changes. The probabilities $p_n$ describe the ordered sizes of blocks of $\Pi_n^0$ at the first time when it differs from $[n]$. Let $D_{1,i}^0$, $2 \leq i \leq n$, be the first time in this homogeneous fragmentation at which 1 and $i$ belong to separate fragments. Let $(\lambda^0(t), t \geq 0)$ be the decreasing process of masses of fragments containing 1. The law of $\lambda^0 = \exp(-\xi)$ is determined by $\Lambda$, as well as that of $(\lambda^0, D_{1,2}^0, \ldots, D_{1,n}^0)$ since

$$\mathbb{P}(D_{1,2}^0 > s_2, \ldots, D_{1,n}^0 > s_n \mid \lambda^0) = \lambda^0(s_2) \cdots \lambda^0(s_n),$$

for all sequences of times $(s_2, \ldots, s_n)$. In particular, knowing $\Lambda$, we know the probabilities $\mathbb{P}(D_{1,2}^0 < \min_{3 \leq i \leq n} D_{1,i}^0) = p_n(1, n-1)$. In the particular case when $n = 3$, this gives $p_3(1, 2)$ and then $p_3(1, 1, 1)$, since $3p_3(1, 2) + p_3(1, 1, 1) = 1$. □

REMARK. It is not hard to see, with a specific example, that in general $\Lambda$ does not characterize the probabilities $p_4(n_1, \ldots, n_k)$, $n_1 + \cdots + n_k = 4$.

PROOF OF PROPOSITION 13. The dislocation measure is determined, up to a scaling constant, by the probabilities $p_n(n_1, \ldots, n_k)$, $\forall n \geq 2$, $\forall (n_1, \ldots, n_k)$ composition of $n$ with $k \geq 2$. The scaling constant is then obtained from $\Lambda$, using (5). The goal here is therefore to check that under the re-rooting assumption, all the probabilities $p_n$ can be recovered from $\Lambda$. Suppose $\Lambda$ is known. We proceed by induction on $n$. For $n=2$, $p_2(1,1)=1$. For $n=3$, the probabilities $p_3$ are known, by Lemma 16. Suppose now that the $p_m$'s are known, $\forall m \leq n-1$. By Lemma 16, $p_n(1, n-1)$ is also known. Then, by Lemma 15, $\forall (n_2, \ldots, n_k)$ composition of $n-2$,

$$p_n(2, n_2, \ldots, n_k) p_2(1,1) = p_n(1, n-1) p_{n-1}(1, n_2, \ldots, n_k),$$

which gives $p_n(2, n_2, \ldots, n_k)$. The probabilities $p_n(n_1, \ldots, n_k)$, with $n_1 \geq 3$, are obtained in the same manner, by induction on $n_1$, thanks to Lemma 15 [note that $p_{n_1}(1, n_1 - 1) \neq 0$, $\forall n_1$]. Therefore, for all compositions $(n_1, \ldots, n_k) \neq (1, \ldots, 1)$, $k \geq 2$, of $n$, we have $p_n(n_1, \ldots, n_k)$, since there is at least one $n_i \neq 1$ and, by symmetry, one can suppose it is $n_1$. It remains to get $p_n(1, \ldots, 1)$, which can be done by using equality (32). $\square$

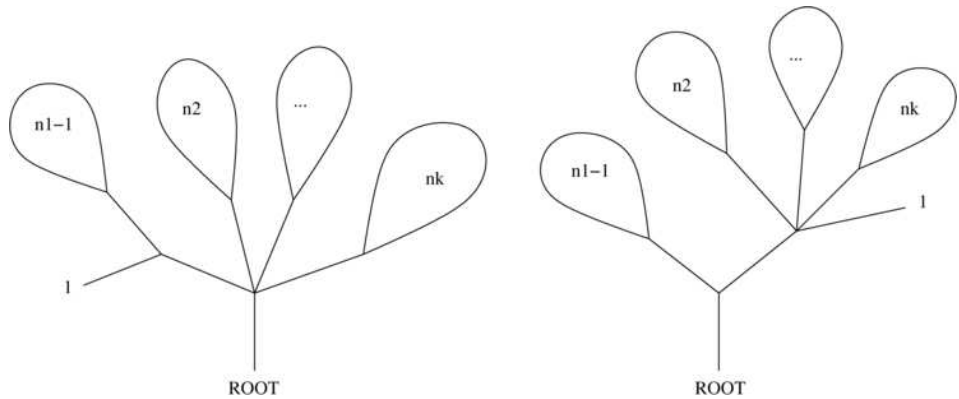

FIG. 5. *By the invariance under re-rooting assumption, these two configurations are equally likely to occur.*

PROOF OF THEOREM 11(ii). By Corollary 14, since the law of the CRT $(\mathcal{T}, \mu)$ is invariant under uniform re-rooting, there exists some $\beta \in (1, 2]$ and some constant $C$ such that $\nu = C\nu_\beta$. It remains to prove that the index of self-similarity is $a = 1/\beta - 1$. Up to now, we only used the combinatorial properties of reduced trees encoded in the dislocation measure $\nu$, and not the further structure of the CRT $(\mathcal{T}, \mu)$ that involves the edge lengths and depends on the scaling parameter $a$. To conclude that $a = 1/\beta - 1$, we must consider edge lengths.

Given the CRT $(\mathcal{T}, \mu)$ rooted at $\rho$ and the leaves $L_1, L_2$, the reduced tree $\mathcal{R}(\mathcal{T}, L_1, L_2)$ can be described by the edge-lengths $D_{1,2}, D_1 - D_{1,2}, D_2 - D_{1,2}$, where $D_{1,2}$ is the separation time of 1 and 2 in $\Pi$ and $D_i$ the first time at which the block containing $i$ is reduced to a singleton, $i = 1, 2$. By invariance under re-rooting, $D_{1,2}$ must have the same law as $D_1 - D_{1,2}$. We already know that this is true for the index $1/\beta - 1$, from Duquesne–Le Gall's Theorem 4.8 in [15].

Using time-changes relating $\Pi$ and its homogeneous counterpart $\Pi^0$ (these time-changes are given specifically in [9]), we have

$$D_{1,2} = \int_0^{D_{1,2}^0} |\Pi_{(1)}^0(t)|^{-a} dt = \int_0^\infty |\Pi_{(1)}^0(t)|^{-a} dt - \int_{D_{1,2}^0}^\infty |\Pi_{(1)}^0(t)|^{-a} dt$$

$$= D_1 - \int_{D_{1,2}^0}^\infty |\Pi_{(1)}^0(t)|^{-a} dt$$

and

(36) $$D_1 - D_{1,2} = \int_{D_{1,2}^0}^\infty |\Pi_{(1)}^0(t)|^{-a} dt,$$

where $D_{1,2}^0$ is the first separation time of 1 and 2 in $\Pi^0$. By the strong Markov property of $\Pi$ (see [9]), $|\Pi_{(1)}^0(t + D_{1,2}^0)|$ has same distribution as $|\Pi_{(1)}^0(D_{1,2}^0)| \cdot |\widetilde{\Pi}_{(1)}^0(t)|$, where $\widetilde{\Pi}^0$ is an independent copy of $\Pi^0$. Therefore,

$$\int_{D_{1,2}^0}^\infty |\Pi_{(1)}^0(t)|^{-a}\,dt = |\Pi_{(1)}^0(D_{1,2}^0)|^{-a} \int_0^\infty |\widetilde{\Pi}_{(1)}^0(t)|^{-a}\,dt = |\Pi_{(1)}^0(D_{1,2}^0)|^{-a} \widetilde{D}_1,$$

where $\widetilde{D}_1$ has same distribution as $D_1$ and is independent of $|\Pi_{(1)}^0(D_{1,2}^0)|^{-a}$. Assuming that $D_{1,2}$ has same distribution as $D_1 - D_{1,2}$ and taking expectations in (36), we obtain

$$\mathbb{E}[|\Pi_{(1)}^0(D_{1,2}^0)|^{-a}]\mathbb{E}[D_1] = \mathbb{E}[D_1 - D_{1,2}]$$
$$= \mathbb{E}[D_{1,2}] = E[D_1](1 - \mathbb{E}[|\Pi_{(1)}^0(D_{1,2}^0)|^{-a}]).$$

For $a < 0$ we may cancel the common factor of $\mathbb{E}[D_1] < \infty$ ($D_1$ is an exponential functional of a subordinator). It remains to notice that the function $f(a) = \mathbb{E}[|\Pi_{(1)}^0(D_{,12}^0)|^{-a}]$ is a strictly monotone function with limit 0 at $-\infty$ and 1 at 0, so that the equation $f(a) = 1 - f(a)$ has a unique solution $a$, which has to be the index $a = 1/\beta - 1$. $\square$

**Acknowledgment.** We thank Grégory Miermont for many interesting discussions.

## REFERENCES

[1] ALDOUS, D. (1991). The continuum random tree. I. *Ann. Probab.* **19** 1–28. MR1085326
[2] ALDOUS, D. (1991). The continuum random tree. II. An overview. In *Stochastic Analysis (Durham, 1990). London Mathematical Society Lecture Note Series* **167** 23–70. Cambridge Univ. Press, Cambridge. MR1166406
[3] ALDOUS, D. (1993). The continuum random tree. III. *Ann. Probab.* **21** 248–289. MR1207226
[4] ALDOUS, D. (1996). Probability distributions on cladograms. In *Random Discrete Structures (Minneapolis, MN, 1993). IMA Vol. Math. Appl.* **76** 1–18. Springer, New York. MR1395604
[5] ALDOUS, D., MIERMONT, G. and PITMAN, J. (2004). Brownian bridge asymptotics for random $p$-mappings. *Electron. J. Probab.* **9** 37–56 (electronic). MR2041828
[6] ALDOUS, D. and PITMAN, J. (1998). Tree-valued Markov chains derived from Galton–Watson processes. *Ann. Inst. H. Poincaré Probab. Statist.* **34** 637–686. MR1641670
[7] BASDEVANT, A.-L. (2006). Ruelle's probability cascades seen as a fragmentation process. *Markov Process. Related Fields* **12** 447–474. MR2246260
[8] BERTOIN, J. (2001). Homogeneous fragmentation processes. *Probab. Theory Related Fields* **121** 301–318. MR1867425
[9] BERTOIN, J. (2002). Self-similar fragmentations. *Ann. Inst. H. Poincaré Probab. Statist.* **38** 319–340. MR1899456

[10] BERTOIN, J. (2006). *Random Fragmentation and Coagulation Processes. Cambridge Studies in Advanced Mathematics* **102**. Cambridge Univ. Press, Cambridge. MR2253162

[11] BERTOIN, J. and ROUAULT, A. (2005). Discretization methods for homogeneous fragmentations. *J. London Math. Soc. (2)* **72** 91–109. MR2145730

[12] DONNELLY, P. and JOYCE, P. (1991). Consistent ordered sampling distributions: Characterization and convergence. *Adv. in Appl. Probab.* **23** 229–258. MR1104078

[13] DUQUESNE, T. and LE GALL, J.-F. (2008). Re-rooting for Lévy trees and snakes. In preparation.

[14] DUQUESNE, T. and LE GALL, J.-F. (2002). Random trees, Lévy processes and spatial branching processes. *Astérisque* **281** 1–147. MR1954248

[15] DUQUESNE, T. and LE GALL, J.-F. (2005). Probabilistic and fractal aspects of Lévy trees. *Probab. Theory Related Fields* **131** 553–603. MR2147221

[16] EVANS, S. N., PITMAN, J. and WINTER, A. (2006). Rayleigh processes, real trees, and root growth with re-grafting. *Probab. Theory Related Fields* **134** 81–126. MR2221786

[17] EVANS, S. N. and WINTER, A. (2006). Subtree prune and regraft: A reversible real tree-valued Markov process. *Ann. Probab.* **34** 918–961. MR2243874

[18] FORD, D. J. (2008). Probabilities on cladograms: Introduction to the alpha model. Preprint. Available at arXiv:math.PR/0511246.

[19] GNEDIN, A. and PITMAN, J. (2005). Regenerative composition structures. *Ann. Probab.* **33** 445–479. MR2122798

[20] GNEDIN, A., PITMAN, J. and YOR, M. (2006). Asymptotic laws for compositions derived from transformed subordinators. *Ann. Probab.* **34** 468–492. MR2223948

[21] GNEDIN, A. V. (1997). The representation of composition structures. *Ann. Probab.* **25** 1437–1450. MR1457625

[22] GREENWOOD, P. and PITMAN, J. (1980). Construction of local time and Poisson point processes from nested arrays. *J. London Math. Soc. (2)* **22** 182–192. MR579823

[23] HAAS, B. (2003). Loss of mass in deterministic and random fragmentations. *Stochastic Process. Appl.* **106** 245–277. MR1989629

[24] HAAS, B. (2007). Fragmentation processes with an initial mass converging to infinity. *J. Theoret. Probab.* **20** 721–758. MR2359053

[25] HAAS, B. and MIERMONT, G. (2004). The genealogy of self-similar fragmentations with negative index as a continuum random tree. *Electron. J. Probab.* **9** 57–97 (electronic). MR2041829

[26] HAAS, B., MIERMONT, G., PITMAN, J. and WINKEL, M. (2008). Continuum tree asymptotics of discrete fragmentations and applications to phylogenetic models. *Ann. Probab.* **36** 1790–1837. MR2440924

[27] KINGMAN, J. F. C. (1978). The representation of partition structures. *J. London Math. Soc. (2)* **18** 374–380. MR509954

[28] LE GALL, J.-F. and WEILL, M. (2006). Conditioned Brownian trees. *Ann. Inst. H. Poincaré Probab. Statist.* **42** 455–489. MR2242956

[29] MIERMONT, G. (2003). Self-similar fragmentations derived from the stable tree. I. Splitting at heights. *Probab. Theory Related Fields* **127** 423–454. MR2018924

[30] MIERMONT, G. (2005). Self-similar fragmentations derived from the stable tree. II. Splitting at nodes. *Probab. Theory Related Fields* **131** 341–375. MR2123249

[31] PITMAN, J. (1999). Coalescents with multiple collisions. *Ann. Probab.* **27** 1870–1902. MR1742892

[32] PITMAN, J. (2006). *Combinatorial Stochastic Processes. Lecture Notes in Math.* **1875**. Springer, Berlin. MR2245368
[33] PITMAN, J. and YOR, M. (1992). Arcsine laws and interval partitions derived from a stable subordinator. *Proc. London Math. Soc. (3)* **65** 326–356. MR1168191
[34] PITMAN, J. and YOR, M. (1997). The two-parameter Poisson–Dirichlet distribution derived from a stable subordinator. *Ann. Probab.* **25** 855–900. MR1434129
[35] SCHROEDER, E. (1870). Vier combinatorische Probleme. *Z. f. Math. Phys.* **15** 361–376.
[36] STANLEY, R. P. (1999). *Enumerative Combinatorics. Vol. 2. Cambridge Studies in Advanced Mathematics* **62**. Cambridge Univ. Press, Cambridge. MR1676282

B. HAAS
CEREMADE
UNIVERSITÉ DE PARIS DAUPHINE
PLACE DU MARÉCHAL
  DE LATTRE DE TASSIGNY 75775 PARIS CEDEX 16
FRANCE
E-MAIL: haas@ceremade.dauphine.fr

J. PITMAN
DEPARTMENT OF STATISTICS
UNIVERSITY OF CALIFORNIA
BERKELEY, CALIFORNIA 94720
USA
E-MAIL: pitman@stat.berkeley.edu

M. WINKEL
DEPARTMENT OF STATISTICS
UNIVERSITY OF OXFORD
1 SOUTH PARKS ROAD
OXFORD OX1 3TG
UNITED KINGDOM
E-MAIL: winkel@stats.ox.ac.uk